\documentclass[11pt,a4paper,leqno,twoside]{amsart}
\usepackage{latexsym,amssymb,amsmath}
\usepackage{epsfig}
\input xy
\xyoption{all}

\newtheorem{lemm}{Lemma}[section]
\newtheorem{theo}[lemm]{Theorem}
\newtheorem{coro}[lemm]{Corollary}
\newtheorem{prop}[lemm]{Proposition}

\def\OP2{\mathbb{OP}^2}
\def\CC{{\mathbb C}}
\def\PP{{\mathbb P}}
\def\QQ{{\mathbb Q}}
\def\ZZ{{\mathbb Z}}
\def\ra{\rightarrow}

\def\cO{\mathcal{O}}
\def\cI{\mathcal{I}}

\def\cH{\mathcal{H}}
\def\cB{\mathcal{B}}

\def\cX{\mathcal{X}}

\def\cE{\mathcal{E}}
\def\lra{\longrightarrow}

\begin{document}

\textheight=22cm

\title{Fano manifolds of degree ten \linebreak and EPW sextics}
\author[A. Iliev]{Atanas Iliev}
\address{Institute of Mathematics,
Bulgarian Academy of Sciences,
Acad. G. Bonchev Str., bl. 8, 
1113 Sofia, Bulgaria}
\email{{\tt ailiev@math.bas.bg}}
\author[L. Manivel]{Laurent Manivel}
\address{Institut Fourier,  
Universit\'e de Grenoble I et CNRS,
BP 74, 38402 Saint-Martin d'H\`eres, France}
\email{{\tt Laurent.Manivel@ujf-grenoble.fr}}

\begin{abstract}
O'Grady showed that certain special sextics in $\PP^5$ called EPW sextics
admit smooth double covers with a holomorphic symplectic structure. We 
propose another perspective on these symplectic manifolds, by showing 
that they can be constructed from the Hilbert schemes of conics on  Fano 
fourfolds of degree ten. As applications, we construct families of Lagrangian 
surfaces in these symplectic fourfolds, and related integrable systems 
whose fibers are intermediate Jacobians.
\end{abstract}

\maketitle

\section{Introduction}

EPW sextics (named after their discoverers, Eisenbud, Popescu and Walter) are some special
hypersurfaces of degree six in $\PP^5$, first introduced in \cite{epw} as examples of Lagrangian 
degeneracy loci. These hypersurfaces are singular in codimension two, but O'Grady realized in 
\cite{EPW1,EPW2,EPW3} that they admit smooth double covers which are irreducible holomorphic symplectic fourfolds. 
In fact, the first examples of such double covers were discovered by Mukai in \cite{mukai2}, who
constructed them as moduli spaces of stable rank two vector bundles on a polarized K3 surface
of genus six. From this point of view, the symplectic structure is induced from the K3 surface. 
It carries over to double covers of EPW sextics by a deformation argument.

The main goal of this paper is to provide another point of view on this symplectic structure. 
Our starting point will be smooth Fano fourfolds $Z$ of index two, obtained by cutting 
the six dimensional Grassmannian
$G(2,5)$, considered in its Pl\"ucker embedding, by a hyperplane and a quadric. Our main observation 
is that the Hodge number $h^{3,1}(Z)$ equals one (Lemma \ref{hodge}). 
By the results of e.g. \cite{km}, a generator 
of $H^{3,1}(Z)$ induces a closed holomorphic two-form on the smooth part of any Hilbert scheme of curves 
on $Z$. We focus on the case of conics. The most technical part
of the paper consists in proving that for $Z$ general, the Hilbert scheme $F_g(Z)$ of conics in $Z$ 
is smooth (Theorem \ref{smoothHs}). It is thus endowed 
with a canonical (up to scalar) global holomorphic two-form. 

Since $F_g(Z)$ has dimension five, it can certainly not be a symplectic variety. However, 
it admits a natural map to a sextic hypersurface $Y_Z^\vee$ in $\PP^5$. We consider the 
Stein factorization 
$$ F_g(Z)\ra \tilde{Y}_Z^\vee\ra Y_Z^\vee.$$
It turns out that $\tilde{Y}_Z^\vee$ is a smooth fourfold, over which $F_g(Z)$ is 
essentially a smooth fibration 
in projective lines. Thus the two-form on $F_g(Z)$
descends to $\tilde{Y}_Z^\vee$. We show that this makes of  $\tilde{Y}_Z^\vee$ a holomorphic symplectic 
fourfold (Theorem \ref{tildesmooth}). Moreover   
the map $\tilde{Y}_Z^\vee\ra Y_Z^\vee$ is a double cover, such that the 
associated involution of $\tilde{Y}_Z^\vee$ is anti-symplectic. This implies that 
$Y_Z^\vee$ is an EPW sextic (Proposition \ref{2EPW}), and that $\tilde{Y}_Z^\vee$ does 
indeed coincide with the double cover constructed by O'Grady (Proposition \ref{sameog}). 

\medskip
Apart from making O'Grady's construction more transparent, at least from our point of view, 
our approach has several interesting consequences. 

First, it shows that double covers of EPW sextics are very close to another classical example 
of symplectic fourfolds, namely the Fano varieties of lines on cubic fourfolds. Indeed, a smooth 
cubic fourfold $Z$ also has $h^{3,1}(Z)=1$, and the symplectic form on its Fano scheme of lines 
$F(Z)$ can be seen as induced from a generator of $H^{3,1}(Z)$, exactly as above. 
Note that a similar line of ideas has been used to explain the existence of a non degenerate
two-form on the symplectic fourfolds in $G(6,10)$ recently discovered in \cite{dv}. 

Second, it sheds some light on the intriguing interplay between the varieties of type 
$Z=G(2,5)\cap Q\cap L$ of different dimensions $N$, where $L$ denotes a linear space of 
dimension $N+4$. If $N=2$, one gets the genus six K3 
surfaces which were, thanks to Mukai's observations, 
at the beginning of that story, but whose associated sextics form 
only a codimension one family in the moduli space of all EPW sextics (see \cite{EPW4}). 
If $N=5$, it is very easy to see that there is 
an EPW sextic attached to $Z$; we explain this in Proposition \ref{dualEPW5}, as a
way to introduce these special sextics. The case $N=3$ was the main theme of investigations 
of \cite{log2} and \cite{dim}; in these studies the surface of conics on $Z$ played a crucial 
role; it is very closely related to the singular locus of the EPW sextic attached to $Z$. 
Finally, for $N=4$ we have seen how to construct an EPW sextic from the family of conics on 
$Z$. In particular, we conclude that for any $N=3, 4$ or $5$, a general EPW sextic is attached 
to a general $Z$, in fact a certain family of such sextics. For sure there is more to understand 
about this, see section 4.5 for a tentative discussion. 

Third, from the fourfold $Z$ we obtain a rather concrete description of the symplectic form on 
$\tilde{Y}_Z^\vee$ (while in \cite{EPW3} its existence was only guaranteed by a deformation 
argument). This allows us to exhibit certain Lagrangian surfaces in $\tilde{Y}_Z^\vee$, that 
we construct either from threefolds that are hyperplane sections of $Z$ (Proposition \ref{lag1}), 
or fivefolds that contain  $Z$ as a hyperplane section (Proposition \ref{lag2}). 
More, we are able to construct, over the moduli stacks parametrizing 
these families of threefolds (respectively, fivefolds), two integrable systems whose Liouville
tori are the corresponding intermediate Jacobians (Theorems \ref{int1} and \ref{int2}). 
Again, this is strikingly similar to the 
constructions of \cite{im}, of two integrable systems  over the moduli stacks parametrizing 
cubic threefolds (respectively, fivefolds) contained in (respectively, containing) a given cubic fourfold.

\medskip\noindent {\it Acknowledgements}. We thank St\'ephane Druel, Dimitri Markushevich and 
Kieran O'Grady for useful discussions. 

\medskip\noindent {\it Notation}.

$V_5$ is a five-dimensional complex vector space. The Grassmannian $G=G(2,5)=G(2,V_5)$ parametrizes 
two-dimensional vector spaces in $V_5$. 

$Z=G\cap Q\cap H$ is the intersection of $G$, considered in its Pl\"ucker embedding, with a quadric $Q$ 
and a hyperplane $H=\PP V_9$, where $V_9\subset\wedge^2V_5$. 

$I_2(Z)$ denotes the linear system of quadrics containing $Z$. The hyperplane of quadrics 
containing $G$, called Pfaffian quadrics, is $I_2(G)\simeq V_5$. 
The hyperplane of Pfaffian quadrics in the
projectivization $I=\PP(I_2(G))\simeq \PP^5$ is denoted $H_P$. 
In the dual projective space $I^\vee$, it defines
a point $h_P$ called the Pl\"ucker point. 

$Y_Z\subset I$ denotes the closure of the locus of singular non Pfaffian quadrics.  
The projectively dual
hypersurface is $Y_Z^\vee\subset I^\vee$. The variety $\hat{Y}_Z^\vee$ parametrizes pairs $(h,V_4)\in 
I^\vee\times\PP V_5^\vee$ such that a quadric in $h$ cuts $\PP(\wedge^2V_4)\cap H$ along a singular quadric. 

$F_g(G)$ is the Hilbert scheme of conics in $G$, $F(G)$ is the nested Hilbert scheme of pairs $(c,V_4)\in
F_g(G)\times\PP V_5^\vee$ such that $c\subset G(2,V_4)$. 

$F_g(Z)$ is the Hilbert scheme of conics in $Z$, $F(Z)$ its preimage in $F(G)$. 

For $c$ a generic conic in $Z$, there is a unique $V_4$ such that $(c,V_4)\in F(Z)$, 
and we denote $G_c=G(2,V_4)$, $P_c=G(2,V_4)\cap H$ and $S_c=P_c\cap Q$. In the pencil of quadrics containing
$S_c$, the unique quadric containing the plane $\langle c\rangle$ spanned by $c$ is denoted $Q_{c,V_4}$.

This defines maps $F(Z)\ra \hat{Y}_Z^\vee$ and $F_g(Z)\ra Y_Z^\vee$. The varieties $\bar{Y}_Z^\vee$ and 
$\tilde{Y}_Z^\vee$ are defined by the Stein factorizations 
$F(Z)\ra \bar{Y}_Z^\vee\ra \hat{Y}_Z^\vee$ and $F_g(Z)\ra \tilde{Y}_Z^\vee\ra Y_Z^\vee$.

\section{EPW sextics in duality}

\subsection{Quadratic sections of $G(2,5)$}
Let $V_5$ be a five dimensional complex vector space. Denote by $G=G(2,V_5)
\subset \PP (\wedge^2V_5)$ the Grassmannian of planes in $V_5$, considered 
in the
Pl\"ucker embedding. Let $X=G\cap Q$ be a general quadric section: this 
is a Fano fivefold of index three and degree ten. In the sequel, when we will
talk about a Fano manifold of degree ten, this will always mean a variety 
of this type, or possibly a linear section.  
 
Let $I=|I_X(2)|$ denote the linear system of quadrics containing 
$X$. Then $I\simeq\PP^5$ is generated by $Q$ and the hyperplane 
$H_P=|I_G(2)|$ of Pfaffian quadrics. Note that once we have chosen an
isomorphism $\wedge^5 V_5\simeq\CC$, there is a natural 
isomorphism 
$$V_5\simeq I_G(2), \qquad v\mapsto P_v(x)=v\wedge x\wedge x.$$
To be more precise, $H_P\simeq \wedge^4V_5^\vee\simeq V_5\otimes \det (V_5^\vee).$

The Pfaffian quadrics $P_v$ are all of rank six. Therefore, the 
divisor $D_X$ of degree ten  parametrizing singular quadrics in $I$ 
decomposes as $$D_X=4H_P+Y_X,$$ 
for some sextic hypersurface $Y_X\subset I$. 

On the other hand, consider a  hyperplane $V_4$ of $V_5$. Then the 
Pl\"ucker quadrics cut $\PP (\wedge^2V_4)\subset \PP (\wedge^2V_5)$ along 
the same quadric, namely the Grassmannian $G(2,V_4)$. Therefore the 
quadrics in $|I_X(2)|$ cut out a pencil of quadrics in $\PP (\wedge^2V_4)$. 
If $V_4$ is general, the generic quadric in this pencil is smooth, 
and there is a finite number of hyperplanes in $|I_X(2)|$ restricting 
to singular quadrics in $\PP (\wedge^2V_4)$. This condition defines 
a hypersurface $Y_X^\vee\subset I^\vee$. The following statement is 
essentially contained in \cite{EPW3} (see in particular Propositions 
7.1 and 3.1). 

\begin{prop}\label{dualEPW5}
The two hypersurfaces $Y_X\subset I$ and $Y_X^\vee\subset I^\vee$ are 
projectively dual EPW sextics.  
\end{prop}

First we need to recall briefly the definition of an EPW sextic (for more
details see \cite{epw,EPW3}; the version we give here follows \cite{EPW4}, section 3.2). 
One starts with a six-dimensional vector space  
$U_6$. Then $\wedge^3U_6$ is twenty-dimensional and admits a natural 
non degenerate skew-symmetric form (once we have fixed a generator of 
$\wedge^3U_6\simeq\CC$). Let then $A\subset\wedge^3U_6$ be a ten-dimensional
Lagrangian subspace. The associated EPW sextic $Y_A^\vee\subset \PP(U_6^\vee)$ is
defined as
$$Y_A^\vee=\{H\subset U_6, \quad \wedge^3H\cap A\ne 0\}.$$
If $A$ is general enough, then $Y_A^\vee$ is singular exactly along 
$$S_A=\{H\subset U_6, \quad \dim(\wedge^3H\cap A)\ge 2\},$$
which is a smooth surface. 

\proof The quadric $Q$ in $G(2,V_5)$ is defined by a tensor in 
$S^2(\wedge^2V_5)^\vee$ modded out by 
the space of Pfaffian quadrics. We choose a representative $Q_0$
in $S^2(\wedge^2V_5)^\vee$. In particular, the 
choice of $Q_0$ induces a decomposition 
$I=H_P\oplus\CC Q_0$, hence a decomposition 
$$\wedge^3I\simeq \wedge^3H_P\oplus \wedge^2H_P\otimes Q_0.$$
Observe that if we let $D=\det V_5^\vee$, then $H_P\simeq V_5\otimes D$, 
hence $\wedge^2H_P\simeq \wedge^2V_5\otimes D^2$ and $\wedge^3H_P\simeq 
\wedge^3V_5\otimes D^3\simeq \wedge^2V_5^\vee\otimes D^2$. 
We can therefore attach to $Q_0$ the subspace $A(Q_0)$ of $\wedge^3I$
defined as 
$$A(Q_0):= \left\{ (Q_0(x,\bullet)\otimes d^2, x\otimes d^2\otimes Q_0), \quad 
x\in\wedge^2V_5 \right\},$$
where $d$ is some generator of $D$. Then $A(Q_0)$ is a Lagrangian subspace 
of $\wedge^3I$, canonically attached to the point defined by $Q_0$ in 
$I-H_P\simeq \CC^5$. 
Consider the EPW sextic $Y_{A(Q_0)}^\vee\subset I^\vee$. 

\begin{lemm} 
$Y_{A(Q_0)}^\vee\simeq Y_X^\vee$. 
\end{lemm}

\proof We prove that $Y_{A(Q_0)}\supset Y_X^\vee$. Since they are both 
sextic hypersurfaces, this will imply the claim. 

A point of $Y_X^\vee$ is defined by a hyperplane $H\subset I$ parametrizing 
quadrics that are all singular when restricted to $\PP(\wedge^2V_4)$, for 
some hyperplane $V_4\subset V_5$. If $H$ is not the Pfaffian hyperplane
$H_P$, we can define it as the space of quadrics of the form 
$Q_v:=P_v-\lambda(v)Q_0$, for some linear form $\lambda$ on $V_5$. 
By the hypothesis, there exists some non zero $p\in \wedge^2V_4$ such that 
$Q_v(p,q)=0$ for any $q\in \wedge^2V_4$. Generically, this $p$ will
not be contained in the cone over $G(2,V_4)$. Otherwise said, $p$ has rank
four, $p\wedge p\ne 0$, and $V_4$ is defined uniquely by $p$.

Observe that the kernel of $\lambda$ must be $V_4$. Indeed, if 
$\lambda(v)=0$, we get that $P_v(p,p)=v\wedge p\wedge p=0$.
But this implies that $v$ belongs to $V_4$. 

The subspace $\wedge^3H$ of $\wedge^3I$ is generated by the tensors
$P_u\wedge P_v\wedge (P_w-\lambda(w)Q_0)$, for $u,v\in V_4$ and $w\in V_5$. 
We can see it as the graph $\Gamma$ of the map 
$\wedge^3H_P\ra\wedge^2H_P\otimes Q_0$ 
induced by the map $H_P\ra\CC Q_0$ sending $P_v$ to $\lambda(v)Q_0$. 
The Lemma is a consequence of the following assertion. 

\medskip\noindent {\it Claim}. The point $(Q_0(p,\bullet), p\otimes Q_0)$ 
belongs to $\Gamma\cap A(Q_0)$. 

\smallskip This is clearly a point of $A(Q_0)$, so we just need to check
that it belongs to $\Gamma$. Observe that $\Gamma$ contains the points 
$(p\wedge w,\lambda(w)p\otimes Q_0)$, for all $w\in V_5$, so we just need 
to prove that there exists some non zero $w$ such that 
$$R_w(\bullet):=p\wedge w\wedge\bullet - \lambda(w)Q_0(p,\bullet)=0.$$
Here $R_w$ is to be considered as a linear form on $\wedge^2V_5$, 
and we know that it 
vanishes on $\wedge^2V_4$. But the orthogonal to $\wedge^2V_4$ in 
$\wedge^2V_5$ is isomorphic to $V_4^\vee$ (once we have chosen a 
generator of $V_5/V_4$), which means that we can identify $R_w$ with 
a linear form $r(w)$ on $V_4$, depending linearly on $w$. But then
the linear map $r: V_5\ra V_4^\vee$ must have a non trivial kernel, 
and we are done. \qed

\begin{lemm} 
$Y_X^\vee$ is dual to $Y_X$. 
\end{lemm}

\proof Consider a general point of $Y_X$, defined by a non Pfaffian 
singular quadric $Q$, with singular point $p$. Suppose that the 
infinitesimally near quadric $Q+\delta Q$ remains singular at the point 
$p+\delta p$. We may suppose that $\delta Q=P_{\delta v}$ and we 
get the order one condition that 
$$(Q+P_{\delta v})(p+\delta p,\bullet )=Q(\delta p,\bullet )+
P_{\delta v}(p,\bullet )=0.$$
Generically, $p\in\wedge^2V_5$ has rank four, hence belongs to 
$\wedge^2V_4$ for a unique hyperplane $V_4$ of $V_5$. We claim
that the quadrics  $Q+\delta Q$ are all singular at $p$, after 
restriction to $\wedge^2V_4$. That is, we claim that 
$$(Q+\delta Q)(p,q)=0 \qquad \forall q\in \wedge^2V_4.$$
Indeed, $\wedge^4V_4$ is one dimensional and generated by $p\wedge p$, 
hence $p\wedge q=\alpha(q)p\wedge p$ for some linear form $\alpha$
on $\wedge^2V_4$. Then, by the identity above, and the fact that 
$Q(p,\bullet)=0$ since $Q$ is singular at $p$, we get 
$$\begin{array}{lcl}
 (Q+\delta Q)(p,q) & = & \delta Q(p,q) \\
 &= & P_{\delta v}(p,q) \\
& = & \delta v\wedge p\wedge q \\
& =& \alpha(q)\delta v\wedge p\wedge p \\
& =& \alpha(q) P_{\delta v}(p,p) \\
& =& - \alpha(q)Q(\delta p,p) \\
&= & 0.
\end{array}$$
This means that the generic tangent hyperplane to $Y_X$ defines 
a point of $Y_X^\vee$, hence that $Y_X^\vee$ is projectively dual 
to $Y_X$. 

This concludes the proof of the Lemma, and of the Proposition as well. \qed

\subsection{Variants}
Consider now a smooth degree ten variety $X$ of dimension $5-k$, defined 
as the intersection of $G(2,V_5)$ with a quadric and a codimension $k$ 
linear subspace $\PP V_{10-k}$ of $\PP(\wedge^2V_5)$. As before we denote 
by $I\simeq\PP^5$ the linear system of quadrics in $\PP V_{10-k}$ 
containing $X$, and by $H_P$ the Pfaffian hyperplane. Generically the Pfaffian
quadrics have rank $6$, hence corank $4-k$. Hence the 
hypersurface $D_X$ of degree $10-k$, parametrizing  singular quadrics in $I$,
decomposes as 
$$D_X=(4-k)H_P+Y_X,$$
where $Y_X$ is again a sextic hypersurface.

As before, we can also define a hypersurface $Y_X^\vee\subset I^\vee$,
parametrizing the hyperplanes in $I$ made of quadrics whose restrictions
to some $\PP(\wedge^2V_4\cap V_{10-k})$ are all singular, $V_4$ being a hyperplane
in $V_5$. The same proof as for the $k=0$ case yields the following result.  

\begin{prop}\label{sextics}
The two hypersurfaces $Y_X\subset I$ and $Y_X^\vee\subset I^\vee$ are 
projectively dual sextics.  
\end{prop}

\proof The only thing we have to prove is that $Y_X^\vee$ has degree 
six. For this we describe, following \cite{log2}, this hypersurface 
as the image of a degeneracy locus, defined as follows. Consider over
$\PP=\PP V_5^\vee$ the rank two vector bundle $F=\cO_\PP\oplus \cO_\PP(1)$,
and the rank $6-k$ vector bundle $M$ whose fiber over $V_4$ is $\wedge^2V_4
\cap V_{10-k}$. For a generic $V_{10-k}$ this is indeed a vector bundle,
at least for $0\le k\le 2$. Let $\cO_F(-1)$ be the tautological line 
bundle over $\PP(F)$. There is a morphism of vector bundles
$$\eta : \cO_F(-1)\ra S^2M^\vee ,$$
defined by mapping a pair $(z,v)\in \CC\oplus V_5$ to the restriction 
of the quadric $zQ+P_v$ to $M_{V_4}$. Since this restriction does only 
depend on the class of $v$ modulo $V_4$, this mapping factors through $\eta$. 

The first degeneracy locus $\hat{Y}_X^\vee$ of $\eta$, defined by the 
condition that the resulting quadric be singular, is a divisor linearly 
equivalent to 
$$[\hat{Y}_X^\vee]=2c_1(M^\vee)-(6-k)c_1(\cO_F(-1)).$$
One easily computes that $c_1(M^\vee)=3h$, where $h$ denotes the 
hyperplane class of $\PP (V_5)$. 

Observe that there is a natural map from $\PP(F)$ to $I^\vee$. 
Indeed, a point of $\PP(F)$ over some $V_4$ is of the form 
$[\lambda,\phi]$, for $\lambda\in\CC$ and $\phi$ a linear form 
on $V_5$ vanishing on $V_4$. It defines the hyperplane in $I$ 
consisting of quadrics of the form $zQ+P_v$ for $\lambda z+\phi(v)=0$. 

In fact this map $\PP(F)\ra I^\vee$ is just the blow-up of the 
point $h_P$ in $I^\vee$ defined by the Pl\"ucker hyperplane. 
This yields a basis of the Picard group of $\PP(F)$ consisting
of the exceptional divisor $E$, and the pull-back $H$ of the 
hyperplane class of $I^\vee$. A standard computation yields
$c_1(\cO_F(-1))=-E$ and $h=H-E$. Hence
$$[\hat{Y}_X^\vee]=6h+(6-k)E=6H-kE.$$

But the hypersurface $Y_X^\vee$ is just the image of $\hat{Y}_X^\vee$ 
in $I^\vee$. Therefore, this formula reads as follows: $Y_X^\vee$ is 
a degree six hypersurface having multiplicity $k$ at $h_P$. \qed

\section{Conics on  Fano fourfolds of degree ten}

\subsection{Conics on $G(2,5)$}

Consider the Hilbert scheme $F_g(G)$ parametrizing conics in $G=G(2,V_5)$. 
In order to study this scheme, we first recall that conics in $G$ can be
partitioned into three different classes, according to the type of their 
supporting plane: 

\begin{enumerate}
\item $\tau$-conics are conics spanning a plane which is not contained 
in $G(2,V_5)$; any smooth $\tau$-conic can be parametrized by $(s,t)\mapsto
(sv_1+tv_2)\wedge (sv_3+tv_4)$ for some linearly independent vectors 
$v_1,v_2,v_3,v_4$ in $V_5$;
\item $\sigma$-conics are conics parametrizing lines passing through a common
point;  any smooth $\sigma$-conic can be parametrized by $(s,t)\mapsto
v_1\wedge (s^2v_2+stv_3+t^2v_4)$ for some linearly independent vectors 
$v_1,v_2,v_3,v_4$ in $V_5$;
\item $\rho$-conics are conics parametrizing lines contained in a common
plane; any smooth $\rho$-conic can be parametrized by $(s,t)\mapsto
(sv_1+tv_2)\wedge (sv_2+tv_3)$ for some linearly independent vectors 
$v_1,v_2,v_3$ in $V_5$.
\end{enumerate}

Each of these three classes of conics is partitioned into three orbits  
of $PGL(V_5)$, consisting of smooth conics, singular but reduced conics, 
and double lines. In particular $F_g(G)$ has exactly nine $PGL(V_5)$-orbits. 
Exactly two are closed: the orbits $F_{\rho}^1$ and $F_{\sigma}^1$ 
parametrizing double-lines of type $\rho$ 
and of type $\sigma$. They are isomorphic, respectively, with the partial
flag varieties $F(2,3,V_5)$ and $F(1,3,4,V_5)$; their dimensions are
$8$ and $9$. The incidence diagram is the following one:

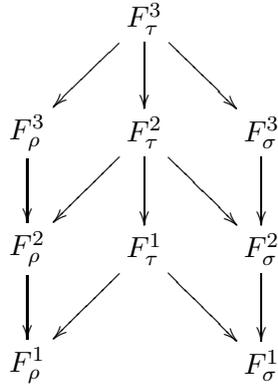
\begin{equation*} 
\xymatrix
{ & F_{\tau}^3\ar[d]\ar[dr]\ar[dl] & \\
 F_{\rho}^3\ar[d] & F_{\tau}^2\ar[d]\ar[dr]\ar[dl] &
F_{\sigma}^3\ar[d] \\
 F_{\rho}^2\ar[d] & F_{\tau}^1\ar[dr]\ar[dl] &
F_{\sigma}^2\ar[d] \\
 F_{\rho}^1 & & F_{\sigma}^1 }
\end{equation*}
 
\medskip

\begin{theo}
The Hilbert scheme $F_g(G)$ of conics in $G=G(2,V_5)$ is irreducible and 
smooth, of dimension $13$. 
\end{theo}

\proof The dimension count is straightforward. 
The singular locus being closed, it is enough to check the smoothness 
at one point of each of the two closed orbits $F_\rho^1$ and $F_\sigma^1$. 
Such a point represents a double-line $\ell$.
Recall (e.g. from \cite{sernesi}) 
that the Zariski tangent space to $F_g(G)$ at the point represented
by $\ell$ is given by 
$$T_{[\ell]}F_g(G)=Hom_{\cO_G}(\cI_\ell,\cO_{\ell}).$$
What we need to check is that this vector space has dimension $13$. 
Since $F_g(G)$ is certainly connected, its smoothness will imply its  
irreducibility.  

\medskip\noindent {\it Double-line of type $\sigma$}. We choose for the 
support of $\ell$ the $\sigma$-plane generated by $v_1\wedge v_2, 
v_1\wedge v_3, v_1\wedge v_4$, for some basis $v_1,\ldots ,v_5$ of 
$V_5$, and we choose in this plane the double line $\ell$ of equations
$$p_{14}^2=0, \quad p_{15}=p_{23}=p_{24}=p_{25}=p_{34}=p_{35}=p_{45}=0,$$
where the $p_{ij}$'s denote the Pl\"ucker coordinates 
on $G(2,V_5)$ associated to our choice of basis. 

We first compute ${\mathcal Hom}_{\cO_G}(\cI_\ell,\cO_{\ell})$ in the 
affine neighborhood of $v_1\wedge v_2$ parametrizing planes which are
transverse to $\langle v_3,v_4,v_5\rangle$. Such a plane has a unique 
basis of the form
\begin{eqnarray*}
 u_1 &= v_1+ x_3v_3+x_4v_4+x_5v_5, \\
 u_2 &= v_2+ y_3v_3+y_4v_4+y_5v_5. 
\end{eqnarray*}
In these coordinates we have $\cI_\ell=\langle y_4^2, y_5, x_3, x_4, x_5\rangle$. An element 
$\phi$ of ${\mathcal Hom}_{\cO_G}(\cI_\ell,\cO_{\ell})$ associates to each of these generators
a section of $\cO_{\ell}$, which can be represented as $p(y_3)+y_4p'(y_3)$ for some polynomials 
$p$ and $p'$. 

We can make the same analysis  in the affine neighborhood of $v_1\wedge v_3$ parametrizing planes which are
transverse to $\langle v_2,v_4,v_5\rangle$. Such a plane has a unique 
basis of the form
\begin{eqnarray*}
 w_1 &= v_1+ z_2v_2+z_4v_4+z_5v_5, \\
 w_3 &= v_3+ t_2v_2+t_4v_4+t_5v_5. 
\end{eqnarray*}
In these coordinates we have $\cI_\ell=\langle t_4^2, t_5, z_2, z_4, z_5\rangle$. An element 
$\psi$ of ${\mathcal Hom}_{\cO_G}(\cI_\ell,\cO_{\ell})$ associates to each of these generators
a section of $\cO_{\ell}$, which can be represented as $q(t_3)+t_4q'(t_3)$ for some polynomials 
$q$ and $q'$. 

Now, we want such a $\psi$ do be defined globally along $\ell$, which means that it must 
extend to a regular morphism over the previous neighborhood. Over $t_2\ne 0$, the formulas
for the change of coordinates are the following:
\begin{eqnarray*}
 x_3=-\frac{z_2}{t_2}, &x_4=z_4-\frac{z_2}{t_2}t_4, &x_5=z_5-\frac{z_2}{t_2}t_5, \\
 y_3= \frac{1}{t_2}, &y_4=\frac{t_4}{t_2}, &y_5=\frac{t_5}{t_2}.
\end{eqnarray*}
Suppose that $\psi$ maps $t_4^2$ to $q(t_2)+t_4q'(t_2)$. Then it maps $y_4^2=t_4^2/t_2^2$
to 
$$t_2^{-2}q(t_2)+t_4t_2^{-2}q'(t_2)=y_3^2q(y_3^{-1})+y_4y_3q'(y_3^{-1}).$$
Therefore $y_3^2q(y_3^{-1})$ and $y_3q'(y_3^{-1})$ must be regular, which means that 
$q$ is at most quadratic and $q'$ is affine. Treating the other conditions similarly 
we check that $\psi$ must be of the following form:
\begin{eqnarray*}
t_4^2 &\mapsto & \psi_1+\psi_2t_2+\psi_3t_2^2+(\psi_4+\psi_5t_2)t_4, \\
t_5 &\mapsto & \psi_6+\psi_7t_2+\psi_8t_4, \\
z_2 &\mapsto & \psi_9+\psi_{10}t_2+\psi_{11}t_4, \\
z_4 &\mapsto & \psi_{12}+\psi_{10}t_4, \\
z_5 &\mapsto & \psi_{13}+\psi_{10}t_4.
\end{eqnarray*}
So there are exactly $13$ free parameters $\psi_1,\ldots ,\psi_{13}$ for $\psi$, as required.

\medskip\noindent {\it Double-line of type $\rho$}. We choose for the 
support of $\ell$ the $\rho$-plane generated by $v_1\wedge v_2, 
v_1\wedge v_3, v_2\wedge v_3$, for some basis $v_1,\ldots ,v_5$ of 
$V_5$, and we choose in this plane the double line $\ell$ of equations
$$p_{23}^2=0, \quad p_{14}=p_{24}=p_{15}=p_{25}=p_{34}=p_{35}=p_{45}=0.$$
 
We compute in the same affine neighborhoods of $v_1\wedge v_2$ and 
$v_1\wedge v_3$. In the latter,  we have $\cI_\ell=\langle z_2^2, z_4, z_5, t_4, t_5\rangle$. An element 
$\psi$ of ${\mathcal Hom}_{\cO_G}(\cI_\ell,\cO_{\ell})$ associates to each of these generators
a section of $\cO_{\ell}$, which can be represented as $q(t_2)+z_2q'(t_2)$ for some polynomials 
$q$ and $q'$. 

A similar analysis as before shows that to be defined globally, such a morphism $\psi$
must be of the following type:
\begin{eqnarray*}
z_2^2 &\mapsto & \psi_1+\psi_2t_2+\psi_3t_2^2+(\psi_4+\psi_5t_2)z_2, \\
t_4 &\mapsto & \psi_6+\psi_7t_2+\psi_8t_4, \\
t_5 &\mapsto & \psi_9+\psi_{10}t_2+\psi_{11}t_4, \\
z_4 &\mapsto & \psi_{12}+\psi_7t_4, \\
z_5 &\mapsto & \psi_{13}+\psi_{10}t_4.
\end{eqnarray*}
Again there are exactly $13$ free parameters $\psi_1,\ldots ,\psi_{13}$ for $\psi$, as required. 
This concludes the proof. \qed

\medskip\noindent {\it Remark}. One can check that $F_g(G)$ is a spherical variety, which means 
that a Borel subgroup of $PGL(V_5)$ acts transitively on some open subset. 
Moreover the Picard number of $F_g(G)$ is three. Indeed, we can consider the nested Hilbert 
scheme $F(G)$ parametrizing pairs $(c,V_4)$, for $c$ a conic in $G$ and $V_4\subset V_5$ a 
hyperplane such that $c$ is contained in the quadric $G(2,V_4)$. One can check that 
the forgetful map $F(G)\rightarrow F_g(G)$ is the blow-up of the codimension two smooth
variety $F_\rho(G)$ parametrizing $\rho$-conics. In particular $F(G)$ contains two 
disjoint divisors $E_\rho(G)$ and $E_\sigma(G)$, the preimages of the subvarieties 
$F_\rho(G)$ and $F_\sigma(G)$ of $F_g(G)$ parametrizing $\rho$ and $\sigma$ conics, respectively. 
These two divisors can themselves be contracted to the variety $S(G)$ parametrizing pairs
$(P,V_4)$, where $P$ is a projective plane inside $\PP(\wedge^2V_4)$. Of course $S(G)$
is a Grassmann bundle over $\PP V_5^\vee$. The condition that 
$P$ be contained inside $G(2,V_4)$ defines two subvarieties $S_\rho(G)$ and $S_\sigma(G)$
(according to the type of $P$), both of codimension six. Blowing-up $S(G)$ over their 
union gives $F(G)$. This is summarized by the following diagram:

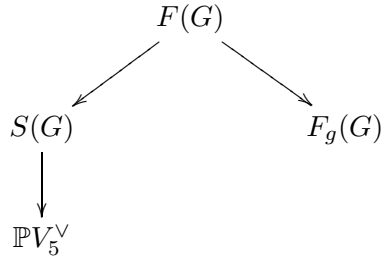
\begin{equation*} 
\xymatrix
{ & F(G)\ar[dr]\ar[dl] & & \\
 S(G)\ar[d] & & F_g(G)   \\
\PP V_5^\vee & & & }
\end{equation*}
 \smallskip

\subsection{Conics on the general Fano fourfold of degree ten}

Now let $Z=G(2,5)\cap H\cap Q$ be a general Fano fourfold of degree ten.
We denote by $F_g(Z)$ the Hilbert scheme of conics in 
$Z$. In this section our main goal is to prove the following statement. 

\begin{theo}\label{smoothHs}
For a general $Z$, the Hilbert scheme $F_g(Z)$ of conics on $Z$ 
is a smooth fivefold. 
\end{theo}

The proof of this result will occupy the rest of the section. 
We will need several auxiliary results, with different techniques 
to handle the three types of conics and their possible singularities. 

\medskip\noindent {\it Reduced conics.}
 
\smallskip We begin with smooth conics. According to its type, the restriction 
to a smooth conic $c$ on $G$, of the dual tautological bundle $T^\vee$, and 
of the quotient bundle $Q$, split as follows (we denote by $\cO_c(1)$ the 
ample generator of the Picard group of $c$, so that $\cO_Z(1)_{|c}=\cO_c(2)$): 

$$\begin{array}{lll}
\mathrm{type} & T_c^\vee & Q_c \\
\tau & \cO_c(1)\oplus \cO_c(1) & \cO_c(1)\oplus \cO_c(1)\oplus \cO_c \\
\sigma & \cO_c(2)\oplus \cO_c & \cO_c(1)\oplus \cO_c(1)\oplus \cO_c \\
\rho  & \cO_c(1)\oplus \cO_c(1) & \cO_c(2)\oplus \cO_c\oplus \cO_c
\end{array}$$

\smallskip
This follows at once from the fact that $T_c^\vee$ and $Q_c$ are globally 
generated and of degree two, and the definitions of the three types of conics. 
This gives the splitting of the tangent bundle $TG=T^\vee\otimes Q$ 
restricted to $c$. We can deduce the splitting type of the normal bundle: 

 $$\begin{array}{ll}
\mathrm{type} & N_{c/G} \\
\tau & \cO_c(2)\oplus \cO_c(2)\oplus\cO_c(2)\oplus \cO_c(1)\oplus \cO_c(1) \\
\sigma & \cO_c(4)\oplus \cO_c(2)\oplus\cO_c(1)\oplus \cO_c(1)\oplus \cO_c \\
\rho  & \cO_c(4)\oplus \cO_c(1)\oplus\cO_c(1)\oplus \cO_c(1)\oplus \cO_c(1)
\end{array}$$

\smallskip
Then we consider the normal exact sequence for the triple 
$c\subset Z\subset G$, 
$$0\ra N_{c/Z}\ra N_{c/G}\stackrel{\theta}{\ra} 
N_{Z/G|c}=\cO_c(2)\oplus \cO_c(4)\ra 0.$$
Our aim is to deduce the possible splitting types of $N_{c/Z}$, and conclude 
that $H^1(N_{c/Z})=0$. This will ensure the smoothness of $F_g(Z)$ at $[c]$.

\begin{lemm}\label{normaltau} 
Let $c\subset Z$ be a smooth $\tau$-conic. Then 
$$N_{c/Z}\simeq \cO_c\oplus\cO_c(1)\oplus\cO_c(1).$$
In particular $F_g(Z)$ is smooth at $[c]$. 
\end{lemm}

 \proof 
There exists a unique hyperplane $V_4\subset V_5$ such that 
$c\subset G_c:=G(2,V_4)$. Moreover $c$ is a linear section of $G_c$, so that 
$N_{c/G_c}=\cO_c(2)\oplus\cO_c(2)\oplus\cO_c(2)$, while 
$N_{G_c/G}=T^\vee_{|G_c}$ (since $G_c$ is the zero locus of a section 
of $T^\vee$ on $G$) and $N_{G_c/G|c}=\cO_c(1)\oplus\cO_c(1)$. 
The normal exact sequence of the triple $c\subset G_c\subset G$ is split. 

Now, $c$ being contained in $Z=G\cap Q\cap H$, it must be 
contained in the quartic surface $S_c=G_c\cap Q\cap H$. 
We get the exact sequence
$$0\ra N_{c/S_c}\ra  N_{c/Z}\ra N_{S_c/Z|c}=N_{G_c/G|c}
=\cO_c(1)\oplus\cO_c(1)\ra 0.$$
But $\omega_{S_c}=\cO_{S_c}(-1)$, hence $\omega_{S_c|c}\simeq \omega_c$. 
Therefore $N_{c/S_c}\simeq\cO_c$ and the exact sequence above must be split. 
This implies the lemma. \qed

\medskip 
Now suppose that $c$ is a $\rho$-conic. Consider in the normal exact 
sequence for the triple $c\subset Z\subset G$, the component $\theta_{44}:
\cO_c(4)\ra \cO_c(4)$ of the morphism $\theta$. We say that $c$ is special
if $\theta_{44}=0$. 

\begin{lemm}
Let $c$ be a non special smooth $\rho$-conic in $Z$. Then 
$$N_{c/Z}\simeq \cO_c\oplus\cO_c(1)\oplus\cO_c(1).$$
In particular $F(Z)$ is smooth at $[c]$. 
\end{lemm}

\proof Since $\theta_{44}\ne 0$, it is an isomorphism and we get an exact 
sequence
$$0\ra N_{c/Z}\ra \cO_c(1)^{\oplus 4}\ra \cO_c(2)\ra 0.$$
In particular $N_{c/Z}^\vee(1)$ is generated by global sections. 
This implies that $N_{c/Z}=\cO_c(n_1)\oplus \cO_c(n_2)\oplus \cO_c(n_3)$
with $n_1, n_2, n_3\le 1$ and $n_1+n_2+n_3=2$. The only possibility is 
that, up to permutation, $n_1=n_2=1$ and  $n_3=0$. \qed

\begin{lemm}
A general $Z$ contains no special $\rho$-conic.
\end{lemm}

\proof One readily checks, with the previous notations,  
that $\theta_{44}=0$ if and only if the quadric 
$Q$ contains the plane spanned by  the $\rho$-conic $c$. 
So, that plane must be contained in $Z$, which is not possible 
for a general $Z$, because of the next easy lemma. \qed
 
\begin{lemm}\label{noplane}
If $Z$ is general, it does not contain any plane.
\end{lemm}

\proof There are two families of planes on $G(2,V_5)$: $\rho$-planes 
of the form $\PP(\wedge^2V_3)$, for $V_3\subset V_5$, and $\sigma$-planes 
of the form $\PP(V_1\wedge V_4)$, for $V_1\subset V_4\subset V_5$. 
The family of $\rho$-planes is parametrized by $G(3,V_5)$, hence 
six-dimensional. The family of $\sigma$-planes is parametrized by the partial
flag variety $F(1,4,V_5)$, hence seven-dimensional. 

Containing a projective plane imposes three conditions on hyperplanes, and 
six conditions on quadrics, hence nine conditions on $Z$. Since nine is bigger 
that seven, we are done. \qed 

\medskip
We can analyze the case of $\sigma$-conics in a similar way: we can define a 
$\sigma$-conic $c$ in $Z$ to be special if  $\theta_{44}=0$. As for  
$\rho$-conic, this implies that the plane spanned by $c$ is contained 
in $Z$, which is not possible for a general $Z$. 

For a non-special $\sigma$-conic $c$, we get an exact sequence
$$0\ra N_{c/Z}\ra \cO_c\oplus \cO_c(1)^{\oplus 2}\oplus \cO_c(2)
\stackrel{\tau}{\ra} \cO_c(2)\ra 0.$$
Again we have two cases, according to the vanishing of the component 
$\tau_{22}: \cO_c(2)\ra \cO_c(2)$. We say that $c$ is of the first kind if 
$\tau_{22}\ne 0$, and of the second kind otherwise. In the latter case, 
$N_{c/Z}=\cO_c(2)\oplus N$, where $N$ fits into an exact sequence
$$0\ra N \ra \cO_c\oplus \cO_c(1)^{\oplus 2}\ra \cO_c(2)\ra 0.$$
This rank two bundle $N$ has degree zero and $N^\vee(1)$ is generated by 
global sections, which leaves only two possibilities: $N=\cO_c\oplus \cO_c$ or
$N=\cO_c(-1)\oplus \cO_c(1)$. We have proved:

\begin{lemm}
Let $c$ be a non special smooth $\sigma$-conic in $Z$. 

If $c$ is of the
first kind,  
$$N_{c/Z}\simeq \cO_c\oplus\cO_c(1)\oplus\cO_c(1).$$

If $c$ is of the second kind, 
$$\begin{array}{ll}
 & N_{c/Z}\simeq \cO_c(2)\oplus\cO_c\oplus\cO_c \\
or & N_{c/Z}\simeq \cO_c(2)\oplus\cO_c(1)\oplus\cO_c(-1).
\end{array}$$

In any case $F(Z)$ is smooth at $[c]$. 
\end{lemm}

This analysis can be extended, with the same conclusions regarding the 
smoothness of $F(Z)$, to reduced singular conics. This was done in 
\cite{IM} in a similar case.  We prefer to present a detailed treatment
of the case of double lines, which requires a different type of arguments. 

\medskip\noindent {\it Double lines.}

\begin{lemm} 
A general $Z$ contains a two-dimensional family of double lines.
This family contains a one dimensional sub-family of double-lines of type 
$\sigma$, and a finite number of double-lines of type $\rho$. 
\end{lemm}

\proof This is just a dimension count. There is an eight-dimensional 
family of lines of $G$, parametrized by the partial flag variety 
$F(1,3,V_5)$. For each line $\ell$, the set of double lines supported
by $\ell$ is parametrized by a projective plane, with a line parametrizing 
double lines of type $\sigma$, and a unique point corresponding to a double line 
of type $\rho$. 

More explicitly, if the line $\ell$ is generated by $e_1\wedge e_2$ 
and  $e_1\wedge e_3$, a double line supported by $\ell$ spans a plane 
$$P=\langle  e_1\wedge e_2,  e_1\wedge e_3, ze_2\wedge e_3+
e_1\wedge f\rangle ,$$
where $f$ is defined up to $\langle  e_1, e_2, e_3\rangle $. 
We can thus parametrize $P$ by the point 
$[z,\bar{f}]\in\PP^2$, where $\bar{f}$ 
denotes the projection of $f$ to $V_5/\langle  e_1, e_2, e_3\rangle $. 
The corresponding double line has type $\sigma$ for $z=0$, and type 
$\rho$ for $\bar{f}=0$. 

In particular we get a ten dimensional family of double lines on the 
Grassmannian $G$. Since containing any of these imposes eight conditions
on $Z$, we are done. \qed

\medskip Now suppose that $\ell$ be a double line in $Z$. Denote by 
$\cI_{\ell,Z}\subset\cO_Z$ its ideal sheaf in $Z$, and by $\cI_{\ell,G}
\subset\cO_G$ its ideal sheaf in $G$. The restriction map $\cI_{\ell,G}
\ra \cI_{\ell,Z}$ induces an exact sequence
$$\begin{array}{ll}
0\ra T_{[\ell]}F_g(Z)=Hom(\cI_{\ell,Z},\cO_\ell)\ra \\
\hspace*{3cm}\ra  
T_{[\ell]}F_g(G)=Hom(\cI_{\ell,G},\cO_\ell)\stackrel{\phi}{\ra}
Hom(\cI_{Z,G},\cO_\ell),
\end{array}$$
where $\cI_{Z,G}$ denotes the ideal sheaf of $Z$ in $G$. Since $\ell$ is 
a smooth point of $F_g(G)$, the Hilbert scheme $F_g(Z)$ is smooth at $[\ell]$
if and only if $\phi$ has rank eight. 

\begin{lemm}\label{double-sing-tau}
Let $\ell$ be a double line of type $\tau$ in $G$. In the variety 
parametrizing the Fano fourfolds $Z$ containing $\ell$, the subvariety 
parametrizing those $Z$ for which $\ell$ is a singular point of $F_g(Z)$,
has codimension at least three.
\end{lemm}

\proof The proof is rather computational, see the Appendix. \qed 

\medskip 
Double lines of type $\sigma$ or $\rho$ can be treated similarly. In fact they 
are easier to handle, since it is enough to show that for such double-lines, 
defining a singular point of $F_g(Z)$ impose at least two, resp.  one, 
conditions  on $Z$. 

This concludes the proof of Theorem 4.3. \qed

\medskip\noindent {\it Remarks.}

1. The variety $F_g^\rho(Z)$ of $\rho$-conics in a general 
$Z=G\cap Q\cap H$  can be analyzed as follows. Since $Z$ contains no plane,
a $\rho$-conic in $Z$ must be the trace of $Q$ over a $\rho$-plane 
of $G$ contained in $H$. Recall that $G\cap H$ can be interpreted as 
an isotropic Grassmannian $IG(2,V_5)$, with respect to a maximal rank 
two-form $\omega$ on $V_5$. All $\rho$-planes in $IG(2,V_5)$ are of 
the form $\PP (V_3)$ for $V_3\subset V_5$ containing the kernel $W_1$ 
of $\omega$. Taking the quotient by this kernel, this identifies the 
variety of $\rho$-planes in $IG(2,V_5)$, 
with a Lagrangian Grassmannian $LG(2,V_4)$, which is nothing else than 
a smooth three-dimensional quadric $\QQ^3$. Hence 
$$F_g^\rho(Z)\simeq \QQ^3.$$ 

2. Similarly, a $\sigma$-conic in $Z$ must be the trace of $Q$ over a 
$\sigma$-plane of $G$ contained in $H$. Such a $\sigma$-plane is defined 
by a flag $V_1\subset V_4$, and it is contained in $H$ if and only if 
$V_4\subset V_1^\perp$, where the orthogonality is taken with respect 
to the two-form $\omega$. There are two cases. If $V_1$ does not coincide
with $W_1$, the kernel of $\omega$, then it determines $V_4$ uniquely. 
If $V_1=W_1$, then $V_4$ can be any hyperplane containing it. One 
easily concludes that 
$$F_g^\sigma(Z)\simeq Bl_0\PP^4,$$
the blow-up of $\PP^4$ at one point. 

\section{A two-form on the Hilbert scheme of conics}

Let $Z=G(2,V_5)\cap Q\cap H$ be a general smooth Fano fourfold of degree ten
and index two. 

\subsection{The Hodge numbers of $Z$}

\begin{lemm}\label{hodge}
The Hodge diamond of $Z$ is the following:

 $$\begin{array}{ccccccccc}
 &&&& 1 &&&& \\
 &&& 0 && 0 &&& \\
 &&0&& 1 &&0 && \\
 &0&& 0 && 0 &&0& \\
0 &&1&& 22 &&1 &&0 \\
 &0&& 0 && 0 &&0& \\
 &&0&& 1 &&0 && \\
 &&& 0 && 0 &&& \\
 &&&& 1 &&&& 
\end{array}$$
\end{lemm}

\proof We write $Z=X\cap Q$ and $X=G(2,V_5)\cap H$, with $H=\PP V_9$. 
In order to compute $h^{3,1}(Z)=h^1(Z,TZ(-2))$, we use the normal sequence
$$0\ra TZ\ra TX_{|Z}\ra \cO_Z(2)\ra 0.$$
The claim that $h^{3,1}(Z)=1$ follows from the fact that $TX(-2)_{|Z}$
has no cohomology in degree zero and one, which itself follows from
the fact that $TX(-2)$ and $TX(-4)$ have no cohomology in degree zero 
and one, and one and two, respectively. But $X$ is a linear section 
of $G$, and by Bott's theorem $TG(-k)$ is acyclic for $1\le k\le 4$. 
This implies that $TG(-k)_{|X}$ is acyclic for $1\le k\le 3$,
and then that $TX(-k)$ is acyclic for $k=2,3$. Finally,  $TG(-5)\simeq
\Omega^5_G$ has non-zero cohomology in degree five only, so $TG(-4)_{|X}$, 
and $TX(-4)$ a fortiori, have no cohomology in degree less than four. 

Now we compute $h^{2,2}(Z)=h^2(Z,\Omega_Z^2)=\chi(Z,\Omega_Z^2)$. 
Observe that the conormal exact sequence of the inclusion $Z\subset G$ 
induces a filtration of $\Omega_{G|Z}^2$ with successive quotients 
$\Omega_Z^2$, $\Omega_Z^1(-1)\oplus \Omega_Z^1(-2)$ and $\cO_Z(-3)$,
so
$$\chi(\Omega_Z^2)=\chi(\Omega_{G|Z}^2)-\chi(\Omega_Z^1(-1))
-\chi(\Omega_Z^1(-2))-\chi(\cO_Z(-3)).$$
Using the Koszul exact sequence we get that 
$$\chi(\Omega_{G|Z}^2)
=\chi(\Omega_G^2)-\chi(\Omega_G^2(-1))-\chi(\Omega_G^2(-2))+
\chi(\Omega_G^2(-3)).$$ 
Bott's theorem yields $\chi(\Omega_G^2)=2$, 
$\chi(\Omega_G^2(-1))=\chi(\Omega_G^2(-2))=0$, and $\chi(\Omega_G^2(-
3))=-5$, hence  $\chi(\Omega_{G|Z}^2)=-3$. Computing the other terms 
similarly we get $\chi(\Omega_Z^2)=22$.  \qed
\medskip

\subsection{The induced form on $F_g(Z)$}

Since $h^{1,3}(Z)=1$, there is a canonical (up to constant) holomorphic 
two-form induced on $F_g(Z)$. At a point defined 
by a smooth conic $c\subset Z$, this two-form can be defined on 
$T_{[c]}F_g(Z)=H^0(N_{c/Z})$ as follows (see \cite{km}). Choose a generator $\sigma$ of 
$H^1(Z,\Omega^3_Z)=H^1(Z,TZ(-2))$. Then consider the composition 
\begin{eqnarray*}
\phi_\sigma : \; \wedge^2H^0(N_{c/Z})\ra H^0(\wedge^2N_{c/Z})=H^0(N^\vee_{c/Z}(2))
\hspace*{2cm}  \\
 { }\hspace*{2cm}\stackrel{\sigma}{\ra} H^1(TZ\otimes N^\vee_{c/Z}(-2))
\ra H^1(\omega_c)=\CC.
\end{eqnarray*}
For the last arrow we used the natural quotient map $TZ_{|c}\ra N_{c/Z}$. Note that
rather than using this map, we can proceed as follows. If $X=G\cap H$, recall
form the proof of Lemma \ref{hodge} 
that a generator of $H^1(Z,\Omega^3_Z)=H^1(Z,TZ(-2))$ is given by 
the extension class of the normal exact sequence 
$$0\ra TZ\ra TX_{|Z}\ra\cO_X(2)\ra 0.$$
On the conic $c$, after dualizing, twisting by $\cO_X(1)_{|c}=\cO_c(2)$ 
and passing 
to the normal bundles of $c$ in $Z$ and $X$, this induces an extension 
\begin{equation}\label{extension}
0\ra \omega_c\ra N_{c/X}^\vee(2)\ra N_{c/Z}^\vee(2)\ra 0.
\end{equation}
We can use directly this extension to produce the map 
$$H^0(N^\vee_{c/Z}(2))\ra H^1(\omega_c)=\CC$$
which defines the two-form at $[c]$, at least up to constant. 

\medskip
Recall that the $\tau$-conic is contained in a unique sub-Grassmannian 
$G(2,V_4)$ of $G(2,V_5)$, and that we denoted by $S_c$ the quartic 
surface $G(2,V_4)\cap H\cap Q$. 

\begin{prop}
Let $c$ be a smooth $\tau$-conic in $Z$. Then the line 
$$H^0(N_{c/S_c})\subset H^0(N_{c/Z})=T_{[c]}F_g(Z)$$ 
is contained in the kernel of $\phi_\sigma$. 
\end{prop}

\proof 
This means that 
the composition of maps above vanishes when restricted to $H^0(N_{c/S_c})
\wedge H^0(N_{c/Z})\subset \wedge^2H^0(N_{c/Z})$. 
Consider the commutative diagram
$$\begin{array}{ccccc}
\wedge^2H^0(N_{c/Z}) & \ra & H^0(\wedge^2N_{c/Z}) & = & H^0(N^\vee_{c/Z}(2))
\\
\uparrow &&\uparrow &&\uparrow 
\\
H^0(N_{c/S_c})\wedge H^0(N_{c/Z}) & \ra & H^0(N_{c/S_c}\wedge N_{c/Z})
& = & H^0(N^\vee_{S_c/Z|c}(2)). 
\end{array}$$
Since $N^\vee_{S_c/Z|c}(2)$ is the restriction to $c$ of the vector bundle
$N^\vee_{S_c/Z}(1)$ on $S_c$, and we can compute the remaining arrows on 
$S_c$ before restricting to $c$. In other words, we can factor through the 
maps
$$
H^0(N^\vee_{S_c/Z}(1)) \stackrel{\cup\sigma}{\ra} 
H^1(TZ\otimes N^\vee_{S_c/Z}(-1))
\ra H^1(\cO_{S_c}(-1)) \ra H^1(\cO_{c}(-2))=\CC.
$$
And the result is clearly zero, since $H^1(\cO_{S_c}(-1))=0$. 
Indeed, this follows
from the Kodaira vanishing theorem when the quartic surface 
$S_c$ is smooth. By continuity, the same conclusion continues to hold
when $S_c$ is singular. \qed

\begin{prop}\label{rank4}
Let $c$ be a generic conic in $Z$. Then $\phi_\sigma$ has rank four 
at the corresponding point $[c]$ of $F_g(Z)$.
\end{prop}

\proof
We denote by  $P_c$ the three-dimensional quadric $G(2,V_4)\cap H$. 
We have $S_c=P_c\cap Q$, and an induced diagram
$$\begin{array}{ccccccc}
 0\ra & \omega_c & \ra & N^\vee_{c/X}(2)& \ra & N^\vee_{c/Z}(2) & \ra 0 \\
 & || & & \downarrow & & \downarrow & \\
 0\ra & \omega_c & \ra & N^\vee_{c/P_c}(2) &\ra  & N^\vee_{c/S_c}(2) & \ra 0.
\end{array}$$
Recall that $N_{c/S_c}\simeq\cO_c$, so that the previous exact sequence 
induces a coboundary map 
$$\kappa_c : H^0(\cO_c(2))\ra H^1(\omega_c)=\CC,$$
which we can consider as a quadratic form on $H^0(\cO_c(1))$. 

\begin{lemm}
The skew-symmetric form $\phi_\sigma$ has rank four at $[c]$ if the quadratic 
form $\kappa_c$ is non degenerate.
\end{lemm}

\proof First recall $N_{S_c/Z}=T^\vee_{S_c}$. 
By the previous proposition, $\phi_\sigma$ factors as 
$$\wedge^2H^0(N_{c/Z})\ra H^0(\wedge^2N_{c/Z})\ra H^0(\wedge^2N_{S_c/Z|c})=
H^0(\cO_c(2))\stackrel{\kappa_c}{\ra}H^1(\omega_c)=\CC.$$
This should be interpreted as follows.
We may suppose that $c$ is a smooth $\tau$-conic, in which case we know 
that $N_{c/Z}=\cO_c\oplus \cO_c(1)\oplus \cO_c(1)$ by Lemma \ref{normaltau}. Hence 
$H^0(N_{c/Z})=\CC\oplus A\oplus A$ if $A=H^0(\cO_c(1))$. 
In this decomposition, the fact that $\phi_\sigma$ factors as we have 
seen means that its matrix is of the form
$$\begin{pmatrix}
 0 & 0 & 0 \\ 0 & 0 & \kappa_c \\
 0 & -\kappa_c & 0 
\end{pmatrix}$$
It is thus clear that $\phi_\sigma$ has rank four if (and only if) 
$\kappa_c$ has rank two. \qed

\medskip What remains to be proved is that, generically $\kappa_c$ is 
non degenerate. For this we can focus on the following situation:
we have a quartic surface $S$ in $\PP^4$ which is a general intersection 
of two quadrics $Q,Q'$ and $c$ is the general conic in $S$. 
We must prove that the exact sequence 
$$0\ra \omega_c\ra N_{c/Q}^\vee(2)\ra N_{c/S}^\vee(2)=\cO_c(2)\ra 0$$
induces a non degenerate quadratic form $\kappa_c$ on $H^0(\cO_c(1))$.

This can be seen as follows. We may suppose that $Q'$ contains the plane
$\langle c\rangle$
spanned by $c$, whose linear equations are, say, $x_3=x_4=0$. This means
that $Q'$ has an equation of the form $x_3m_3+x_4m_4=0$, for some
linear forms $m_3,m_4$. Restricted to $\langle c\rangle$, these linear 
forms define two global sections $q_3,q_4$ of $\cO_c(2)$, and the 
linear form $\kappa_c$ is just the projection map 
$$\kappa_c : H^0(\cO_c(2))\ra H^0(\cO_c(2))/\langle q_3,q_4\rangle \simeq\CC .$$
In other words, $\kappa_c$ is polar to the pencil $\langle q_3,q_4\rangle$. It is
thus non degenerate as soon as this pencil has no base point, which is the general
situation.  \qed

\subsection{The dual sextic}

Recall that we denoted by $I$  the linear system of quadrics 
containing $Z$. We have defined in section 2.2 the hypersurface $Y_Z^\vee$
in $I^\vee$  
as follows. Let $\hat{Y}_Z^\vee\subset I^\vee\times\PP(V_5^\vee)$
be the variety 
parametrizing pairs $(h,V_4)$ such that quadrics in $h\subset I^\vee$
cut $\PP(\wedge^2V_4)\cap H$ along singular quadrics. Then  $Y_Z^\vee$ 
is just the image of $\hat{Y}_Z^\vee$ by the first projection. 

Note that generically, for $(h,V_4)$ in $\hat{Y}_Z^\vee$, quadrics of
the hyperplane $h$ will restrict to a corank one quadric in  
$\PP(\wedge^2V_4)\cap H$. Let $S_Z$ denote the locus where the 
corank of the restricted quadric is bigger than one. 

\begin{prop}\label{singhat}
For $Z$ general, the variety $\hat{Y}_Z^\vee$ is an irreducible 
fourfold whose singular locus is exactly $S_Z$. Moreover $S_Z$ is 
a smooth surface, and $\hat{Y}_Z^\vee$ has multiplicity two at any 
point of $S_Z$. 
\end{prop}

\proof 
As in \cite{log2}, and as we have seen in the proof of Proposition 
\ref{sextics}, the variety $\hat{Y}_Z^\vee$ is a degeneracy 
locus for a section of a bundle of quadrics. The conclusion will thus 
follow from a transversality argument: if the section is general enough, 
such a degeneracy locus $\hat{Y}_Z^\vee$ is singular exactly 
along the next degeneracy locus (see \cite{ACGH}, Chapter 2), 
which is precisely $S_Z$. Since the 
next one has, generically, codimension three in $S_Z$, it is in
fact empty, and $S_Z$ must be smooth. Unfortunately, in our section 
we do not deal with a general section, and we will need to check the 
transversality condition explicitly. 

We recall the setting: over $\PP=\PP V_5^\vee$, first consider 
the rank two vector bundle $F=\cO_\PP\oplus \cO_\PP(1)$,
then the rank five vector bundle $M$ whose fiber over $V_4$ is $\wedge^2V_4
\cap V_9$. Let $\cO_F(-1)$ be the tautological line 
bundle over $\PP(F)$. Then we denoted by 
$$\eta : \cO_F(-1)\ra S^2M^\vee $$
the morphism of vector bundles defined by mapping a pair $(z,v)\in 
\CC\oplus V_5$ to the restriction of the quadric $zQ+P_v$ to $\wedge^2V_4
\cap V_9$.

We will only need to consider quadrics of corank one or two, since:

\begin{lemm}\label{nocorank3}
For a general $Z$, the image of $\eta$ does not contain any quadric 
of corank three or more.
\end{lemm}

\proof A straightforward dimension count. \qed

\medskip {\it Corank one}.

Consider a point of $\hat{Y}_Z^\vee$ defined by a corank one quadric. 
We will prove it is a smooth point of $\hat{Y}_Z^\vee$. If this quadric
is not the Pl\"ucker one, we may suppose, up to a change of notation, 
that $Q$ itself cuts $\wedge^2V_4\cap V_9$ along a corank one quadric, 
singular at $\omega_0$. We choose local coordinates on $\PP(F)$ at the corresponding
point as follows. First, we choose a supplement $V_1$ of $V_4$, so that 
any hyperplane in $V_5$ transverse to $V_1$ can be represented as the graph
$V_4(\phi)$ of a morphism $\phi\in Hom(V_4,V_1)$. Then, we represent the 
trace of a hyperplane in $I$ on $\wedge^2V_4(\phi)\cap V_9$ by the restriction
of the quadric $Q+P_v$, for some $v\in V_1$. Our local coordinates will be the 
pair $(\phi,v)$. 

We want to represent the latter quadric $Q+P_v$ on $\wedge^2V_4(\phi)\cap V_9$ by 
an isomorphic quadric $Q_{\phi,v}$ on $\wedge^2V_4\cap V_9$. To do this, we first
observe that $\phi$ induces an isomorphism from $\wedge^2V_4$ to $\wedge^2V_4(\phi)$
sending $\omega$ to $\phi(\omega):=\omega+\phi\rfloor\omega$ (where the contraction map $\phi\rfloor$
maps $v\wedge v'$ to $\phi(v)\wedge v'+v\wedge \phi(v')$). Let $h$ be an equation 
of $H$, and $\Omega$ in $\wedge^2V_4$ be such that $h(\Omega)=1$. If 
$\omega$ belongs to $\wedge^2V_4\cap V_9$, and $\omega'=\omega+t\Omega$, then
$\phi(\omega')$ belongs to $V_9$ when 
$t=-h(\phi\rfloor\omega),$
up to terms of higher order in $(\phi,v)$. Up to such terms, we thus let 
$$\begin{array}{ccl}
Q_{\phi,v}(\omega) &= &(Q+P_v)(\phi(\omega')) \\
 &= & Q(\omega)+2Q(\omega,\phi\rfloor\omega-h(\phi\rfloor\omega)\Omega)+
v\wedge\omega\wedge\omega .
\end{array}$$
Since $Q$ is supposed to have corank one and kernel $\langle\omega_0\rangle$, the function 
$\det(Q_{\phi,v})$ is equal, up to a constant and higher order terms, to 
$Q_{\phi,v}(\omega_0)$, which is therefore a local equation of $\hat{Y}_Z^\vee$.
For this equation to vanish identically at first order, we would need that 
\begin{eqnarray}
v\wedge\omega_0\wedge\omega_0 &=0 &\forall v\in V_1, \\
Q(\omega_0,\phi\rfloor\omega_0-h(\phi\rfloor\omega_0)\Omega) &=0 &\forall\phi\in Hom(V_4,V_1).
\end{eqnarray}
Since $V_1$ is transverse to $V_4$, the first equation implies that $\omega_0\wedge\omega_0=0$, 
hence $\omega_0$ has rank two and defines a point of $G(2,V_4)\cap H$. We will denote  
the corresponding plane by $V_2\subset V_4$. Observe that, when $\phi$ varies,
$\phi\rfloor\omega_0$ describes $V_1\wedge V_2\subset\wedge^2V_5$. Since, by the singularity condition,
$Q(\omega_0,\omega)=0$ for any $\omega\in\wedge^2V_4\cap V_9$,  the second condition means that
the linear form $Q(\omega_0,\omega)$ is proportional to $h$ on $V_2\wedge V_1\oplus \wedge^2V_4$. 
We claim that this implies that $[\omega_0]$ is a singular point of $Z$, thus leading to a contradiction. 
Indeed, the affine tangent space to $G$ at $\omega_0$ is $V_2\wedge V_5$, which is contained in 
$V_2\wedge V_1\oplus \wedge^2V_4$. So the traces on this tangent space, of $H$ and of the tangent space 
to $Q$, would coincide, and the intersection of $G$, $H$ and $Q$ would not be transverse at $\omega_0$. 

Suppose now that the Pl\"ucker quadric itself, $G(2,V_4)\cap H$, is singular at $\omega_0$. 
Choose a generator $v_1$ of $V_1$. Local coordinates on $\PP(F)$ are given by $(t,\phi)$, 
where $\phi\in Hom(V_4,V_1)$ as above and the quadric to consider on $\wedge^2V_4(\phi)\cap V_9$ is
$P_{v_1}+tQ$. As in the previous case we identify $\wedge^2V_4(\phi)\cap V_9$ with 
$\wedge^2V_4\cap V_9$, and the previous quadric on $\wedge^2V_4(\phi)\cap V_9$ with the 
isomorphic one on $\wedge^2V_4\cap V_9$ given by 
$$Q_{t,\phi}(\omega)=v_1\wedge\omega\wedge\omega+tQ(\omega)-2h(\phi\rfloor\omega_0)v_1\wedge\Omega\wedge\omega,$$
up to order one. As above, we want to exclude the possibility that 
$$Q_{t,\phi}(\omega_0)=tQ(\omega_0)-2h(\phi\rfloor\omega_0)v_1\wedge\Omega\wedge\omega_0=0$$
for all $t$ and $\phi$. Since $\Omega\wedge\omega_0\ne 0$ (otherwise $G(2,V_4)$ would be singular
at $[\omega_0]$ !), this would mean that $Q(\omega_0)=0$ and $h(V_2\wedge V_1)=0$ where, as above,
$V_2$ is the plane defined by $\omega_0$. The former condition means that $\omega_0$ defines a point of $Z$. 
The latter one implies that $H$ is tangent to $G$ at $[\omega_0]$. As in the previous case we would therefore
conclude that $Z$ is singular at $[\omega_0]$, a contradiction.

\smallskip {\it Corank two}.

Since $G(2,V_4)$ is smooth, a hyperplane section cannot have corank two.  
We may therefore suppose that $Q$ itself cuts $\wedge^2V_4\cap V_9$ along a corank two quadric, 
singular along the line $\langle\omega_0,\omega_1\rangle$. Considering the 
intersection of this line with the quadric  $G(2,V_4)\cap H$, we may suppose that 
$\omega_0$ and $\omega_1$ have rank two. We choose local coordinates $(v,\phi)$ 
on $\PP(F)$ as above and we consider the same quadric $Q_{v,\phi}$. The required 
transversality condition can be expressed by the condition that the map 
$$(v,\phi)\mapsto \begin{pmatrix} Q_{v,\phi}(\omega_0,\omega_0) & Q_{v,\phi}(\omega_0,\omega_1)\\
Q_{v,\phi}(\omega_0,\omega_1) & Q_{v,\phi}(\omega_1,\omega_1) \end{pmatrix}$$
have rank three. If $\omega_0\wedge\omega_1\ne 0$, the off diagonal term 
$Q_{v,\phi}(\omega_0,\omega_1)$ will contribute by one to the rank, through its terms involving $v$. 
So what we need to avoid is that $Q(\omega_0,\phi\rfloor\omega_0-h(\phi\rfloor\omega_0)\Omega)$
and $Q(\omega_1,\phi\rfloor\omega_1-h(\phi\rfloor\omega_1)\Omega)$ impose linearly dependent 
conditions. Since we have 
supposed that $\omega_0$ and $\omega_1$ represent transverse planes in $V_4$, 
this must be the case, except if  one of these linear forms is zero. This can be excluded,
for a general $Z$, by a simple dimension count: for a given $H$, we have four parameters
for $V_4$, then three for each of $[\omega_0]$ and $[\omega_1]$, which belong to the 
three-dimensional quadric $G(2,V_4)\cap H$. But we get eleven linear conditions on $Q$. 

Of course we also need to consider the degenerate cases for which the line $\langle\omega_0,\omega_1\rangle$
is tangent to $G(2,V_4)\cap H$, or even contained in  $G(2,V_4)\cap H$ (this would mean that 
$\omega_0\wedge\omega_1=0$). A similar dimension count leads in both cases to the same conclusion. 
\qed 

\begin{prop}\label{blowdownhat}
The projection map $\hat{Y}_Z^\vee\ra Y_Z^\vee$ is the blow-up of the 
point of $Y_Z^\vee$ defined by the Pl\"ucker hyperplane.
\end{prop}

\proof 
This follows from the proof of Proposition \ref{sextics}. Indeed, we 
have seen in this proof that  $\hat{Y}_Z^\vee\ra Y_Z^\vee$ is the 
restriction of the blow-up of the point $h_P$ in $I^\vee$ defined by the 
Pl\"ucker hyperplane. Moreover, we are in the case $k=1$ of the 
Proposition, and the proof shows that $Y_Z^\vee$ has multiplicity 
one at this point. Otherwise said, this is a smooth point of $Y_Z^\vee$. 
This is enough to ensure that the 
projection map $\hat{Y}_Z^\vee\ra Y_Z^\vee$ is just the blowing-up of $h_P$. 

Note that the preimage of $h_P$ in $\hat{Y}_Z^\vee$ is easily determined. 
It is simply the hyperplane in $\PP V_5^\vee$, 
defined as the set of those hyperplanes that 
contain the kernel of the (degenerate) two-form defining $H$. 
\qed

\subsection{The symplectic structure}

Consider the nested Hilbert scheme $F(Z)$ parametrizing pairs $(c,V_4)$ 
such that $c$ be a conic in $Z$ and $V_4$ a hyperplane in $V_5$ such that 
the quadric $G(2,V_4)$ contains $c$. Recall that $V_4$ is uniquely determined
by $c$ except if $c$ is a $\rho$-conic, in which case there is a projective 
line of possible $V_4$'s. Otherwise said, the forgetful map 
$$\pi : F(Z)\ra F_g(Z)$$ 
is an isomorphism outside $F_g^\rho(Z)$, and contracts a divisor 
onto $F_g^\rho(Z)$. 

\begin{lemm}
For $Z$ general, $F(Z)$ is smooth.
\end{lemm}

\proof 
The nested Hilbert scheme $F(Z)$ is a  subscheme of $F_g(Z)\times
\PP V_5^\vee$. Its Zariski tangent space at a point $(c,V_4)$ can be 
described by the following exact sequence, where we let $G_c=G(2,V_4)$:   
$$0\ra T_{(c,V_4)}F(Z)\ra H^0(N_{c/Z})\oplus H^0(N_{G_c/G})
\ra H^0(N_{G_c/G|c}).$$
Since we already know that $F_g(Z)$ is smooth at $c$, the smoothness of 
$F(Z)$ at $(c,V_4)$ is equivalent to the surjectivity of the rightmost 
arrow. 

First observe that $N_{G_c/G}=T^\vee_{|G_c}$, so that $H^0(N_{G_c/G})$ is 
simply $V_4^\vee$, and the restriction map $H^0(N_{G_c/G})
\ra H^0(N_{G_c/G|c})$ is already surjective if $c$ is not a $\rho$-conic. 

So suppose that $c$ be a $\rho$-conic, spanning a plane $G(2,V_3)$ with
$V_4\supset V_3$. In this case the restriction map  $H^0(N_{G_c/G})
\ra H^0(N_{G_c/G|c})$ has for image a hyperplane $H_c$, and we need to 
check that the 
image of the other map $H^0(N_{c/Z})\ra H^0(N_{G_c/G|c})$ is not contained 
in $H_c$. 

This condition amounts to the fat that a matrix of size 
$8\times 4$ be of maximal rank. If this matrix is sufficiently general,
this fails to happen in codimension $8-4+1=5>4$. Our claim follows. \qed

\medskip\noindent {\it Remark}. One can show that the map $\pi : F(Z)\ra F_g(Z)$
is simply the blow-up of the codimension two subvariety $F_\rho(Z)$ parametrizing
$\rho$-conics in $Z$. This subvariety is smooth for $Z$ general enough. 

\medskip
Our next goal is to construct a morphism 
$$\alpha : F(Z)\lra \hat{Y}_Z^\vee.$$
For a point $(c,V_4)$ in $F(Z)$, we have two quadratic hypersurfaces inside 
$\PP(\wedge^2V_4)\cap H\simeq \PP^4$: the intersection $P_{V_4}$ of the 
Pl\"ucker quadric  $G(2,V_4)$ with $H$, and the trace $Q_{V_4}$ of the 
quadric $Q$ defining $Z$. The pencil $\langle P_{V_4}, Q_{V_4}
\rangle$ is uniquely defined by $Z$ (and $V_4$). Any quadric in this pencil 
contains the conic $c$, but the generic one does not contain the plane 
$\langle c\rangle$ spanned by $c$, since that plane cannot be contained in $Z$
by Lemma \ref{noplane}. 
This implies that there is a unique quadric $Q_{c,V_4}\in
\langle P_{V_4}, Q_{V_4}\rangle$ containing $\langle c\rangle$. Moreover 
this quadric must be singular, since a smooth three dimensional quadric 
does not contain any plane. Therefore $Q_{c,V_4}$ defines a point 
in $\hat{Y}_Z^\vee$: this is $\alpha(c,V_4)$. 

\begin{prop}
Let $y$ be a point of $\hat{Y}_Z^\vee$. 
\begin{enumerate}
\item If $y\in S_Z$, the set-theoretical fiber $\alpha^{-1}(y)$ is a projective line.
\item If $y\notin S_Z$, the fiber $\alpha^{-1}(y)$ is a disjoint union
of two projective lines.
\end{enumerate}
\end{prop}

\proof 
Suppose that $y=\alpha(c,V_4)$ does not belong to $S_Z$. This means that 
the quadric $Q_{c,V_4}$ has rank four; otherwise said, it is a cone over
a smooth quadratic surface. The projective planes $L$ in $Q_{c,V_4}$ are then
cones over the lines in this surface, and are parametrized by two projective 
lines. Any such plane $L$, cut out with any other quadric in the pencil 
$\langle P_{V_4}, Q_{V_4}\rangle$, gives a conic $c(L)$ such that 
$Q_{c(L),V_4}=Q_{c,V_4}$, hence $\alpha(c(L),V_4)=\alpha(c,V_4)=y$. 
This implies that $\alpha^{-1}(y)$ is the disjoint union
of these two projective lines.

If $y=\alpha(c,V_4)$ does belong to $S_Z$, the quadric $Q_{c,V_4}$ has 
rank three; it is a double cone over a smooth conic. The projective planes 
$L$ in $Q_{c,V_4}$ are then parametrized by that single conic. We can make 
with these planes the same construction as above, but we end up with a 
single projective line parametrizing  $\alpha^{-1}(y)$. \qed

\medskip\noindent {\it Remark}. Observe what happens over conics not of type
$\tau$. First recall that $F^\sigma(Z)\simeq F_g^\sigma(Z)\simeq Bl_0\PP^4$, 
the blow-up of $\PP^4$ at one point. This blow-up is a $\PP^1$-fibration 
of $\PP^3$, and the restriction of $\alpha$ to $F^\sigma(Z)$ coincides
with this fibration. Second, recall that $F_g^\sigma(Z)\simeq \QQ^3$, 
coincides with the isotropic Grassmannian $IG(3,V_5)$. Its preimage 
$F^\sigma(Z)$ in $F(Z)$ is the variety $IF(3,4,V_5)$ of flags $V_3\subset
V_4$ with $V_3$ isotropic (recall that this implies that $V_3$ contains
$W_1$, the kernel of the two-form $\omega$ on $V_5$ defining $H$). 
The restriction of the map $\alpha$ to $F^\sigma(Z)$ simply forgets 
$V_3$. In particular its image is the space of hyperplanes in $V_5$ 
containing $W_1$, hence a copy of $\PP^3$. 

\begin{prop}
Any fiber of $\alpha$ is a smooth curve in $F(Z)$.
\end{prop}

\proof Let $(c,V_4)$ be a point of $F(Z)$. We want to prove that the corresponding 
fiber $F_c$ of $\alpha$ is smooth at that point. Set-theoretically, we have seen that 
the plane $\langle c\rangle$ is contained in a unique  (singular) quadric
$Q_{c,V_4}$ of the pencil of quadrics obtained by restricting $I$ to $\PP(\wedge^2V_4)\cap H$. This plane 
$\langle c\rangle$ varies in a family of planes in $Q_{c,V_4}$ parametrized by a projective line, and we 
get a map $\PP^1\ra F_c$, which we shall prove to be a local isomorphism.  

Observe that the tangent space to $F_c\subset F(Z)$ is $H^0(N_{c/S})\subset H^0(N_{c/Z})$,
where $S$ is the quartic surface cut out by $Z$ on $\PP(\wedge^2V_4)\cap H$. This surface 
is the intersection $S=Q_0\cap Q_{c,V_4}$ of two quadrics. We can choose linear coordinates $x_0,\ldots ,x_4$ 
such that $\langle c\rangle$ be defined by $x_2=x_3=0$, and write  
\begin{eqnarray*}
Q_0 &= &x_3\ell_3+x_4\ell_4+q(x_0,x_1,x_2), \\
Q_{c,V_4} &= &x_3m_3+x_4m_4, 
\end{eqnarray*}
for some linear forms $\ell_3,\ell_4,m_3,m_4$. Note that  $q(x_0,x_1,x_2)$ (which is 
non zero since $Z$ contains no plane) is an equation of $c$ in $\langle c\rangle$. 

Now we can see very explicitly that the map $T_{\langle c\rangle}\PP^1\ra H^0(N_{c/S})$ is non zero,
which will prove our claim.  Indeed, an infinitesimal deformation of $\langle c\rangle$ in $Q_{c,V_4}$ is simply
obtained by $\epsilon\mapsto \langle c\rangle(\epsilon)$, the plane defined by the two equations
$x_3+\epsilon m_4=x_4-\epsilon m_3=0$. It is mapped to a global element $\theta$ of 
$H^0(N_{c/S})=Hom_{\cO_S}(\cI_c,\cO_c)$ defined by 
\begin{eqnarray*}
 x_3 &\mapsto &m_4, \\
 x_4 &\mapsto &-m_3 \\
 q &\mapsto &m_3\ell_4-m_4\ell_3.
\end{eqnarray*}
Indeed, $\cI_c$ is generated by $x_3, x_4$ and $q$ at any point (strictly speaking, to make sense 
of this we need to divide them by some linear, respectively quadratic form not vanishing at the 
point considered), and although $\ell_3,\ell_4,m_3,m_4$ are not uniquely defined, $m_4$, $m_3$ and 
$m_3\ell_4-m_4\ell_3$ are uniquely defined when restricted to $c$. 

There just remains to check that $\theta$ cannot be zero. This would mean that $m_3$ and $m_4$
vanish identically on $c$, hence that they are linear combinations of $x_3$ and $x_4$. But then
$Q_{c,V_4}$ would have rank at most two, and by Lemma \ref{nocorank3}, this is not possible for a general $Z$.
\qed

\medskip
Consider the Stein factorization of $\alpha$:
$$F(Z)\stackrel{\beta}{\lra} \bar{Y}_Z^\vee
\stackrel{\gamma}{\lra} \hat{Y}_Z^\vee.$$
By the previous proposition, $\gamma$ has degree two, and ramifies precisely
over $S_Z$. By the previous Proposition 
the reduced fibers of $\beta$ are smooth projective lines, and 
in fact $\beta$ is a $\PP^1$-bundle, since an application of \cite[Theorem 4.1]{AW} yields:
  
\begin{prop}
The variety $\bar{Y}_Z^\vee$ is smooth.
\end{prop}


Note that conics in $Z$ which are not $\tau$ conics are sent to the 
Pl\"ucker hyperplane in $Y_Z^\vee$. 
In particular the map $F(Z)\ra Y_Z^\vee$ factorizes through $F_g(Z)$. 
Taking the Stein factorization of the induced map $F_g(Z)\ra Y_Z^\vee$, 
we get a commutative diagram
$$\begin{array}{ccccc}
F(Z) & \ra & \bar{Y}_Z^\vee & \ra & \hat{Y}_Z^\vee \\
\downarrow & & \downarrow & & \downarrow \\
F_g(Z) & \ra & \tilde{Y}_Z^\vee & \ra & Y_Z^\vee 
\end{array}$$

\smallskip
A consequence of the previous lemma is that:

\begin{lemm}
The projection map $\bar{Y}_Z^\vee\ra \tilde{Y}_Z^\vee$ is the blow-up of the 
two points of $ \tilde{Y}_Z^\vee$ in the preimage of the Pl\"ucker hyperplane.
In particular $ \tilde{Y}_Z^\vee$ is smooth.
\end{lemm}

Now we can prove the main result of this section:

\begin{theo}\label{tildesmooth}
The variety $\tilde{Y}_Z^\vee$ is a smooth symplectic fourfold.
\end{theo}

\proof
We can use our two-form $\phi_\sigma$ on $F_g(Z)$ and lift it 
to $F(Z)$. Since the fibers of $\beta$ are projective lines, the 
induced two-form on $F(Z)$ descends to a globally defined two-form 
$\Phi_\sigma$ on $\bar{Y}_Z^\vee$, which remains a closed form.
The generic rank of  $\Phi_\sigma$ is four since the generic rank of  
$\phi_\sigma$ is four by Proposition \ref{rank4}. 
Since the projection to $\bar{Y}_Z^\vee$
is birational, we also get a closed two-form 
$\tilde\Phi_\sigma$ on $\tilde{Y}_Z^\vee$, generically non-degenerate. 

But the canonical class of $\tilde{Y}_Z^\vee$ is trivial, implying 
that $\tilde\Phi_\sigma$ is in fact everywhere non-degenerate. Indeed, the 
sextic $Y_Z^\vee$ is smooth in codimension one, hence normal. Its canonical
class is trivial. The map  $\tilde{Y}_Z^\vee\ra Y_Z^\vee$ is finite of 
degree two, ramified on the surface $S_Z$ only, so the canonical class
of  $\tilde{Y}_Z^\vee$ is simply the pull-back of that of $Y_Z^\vee$.
Hence the claim and the theorem. \qed

\subsection{EPW sextics attached to Fano manifolds of different dimensions}
 
Let us elaborate on what we have proved at this point. Consider a general 
variety $X=G\cap Q\cap \PP V_{10-k}$, of dimension $N=5-k$, and the associated sextic hypersurface
$Y_X$. We have recalled in Proposition \ref{dualEPW5} that for $k=0$, $Y_X$ is an EPW sextic. We have proved 
it is also the case for $k=1$. This is also true for $k=3$, in which case 
$X$ is a generic polarized K3 surface of degree ten \cite{mukai1}. Mukai 
showed (\cite{mukai2}, Ex. 5.17) that the natural double cover $\tilde{Y}_X$ of the 
sextic $Y_X$ can be identified with the moduli space of stable rank two vector 
bundles $E$ on $X$ with Chern classes $c_1(E)=\cO_X(1)$ and $c_2(E)=5$. This explains the 
existence of a symplectic structure on $\tilde{Y}_X$, directly inherited from that of $X$. 

\medskip\noindent {\it Remark}. 
O'Grady proved that in the (irreducible) family of EPW sextics, those coming from polarized 
K3 surfaces of degree ten form a codimension one family (\cite{EPW4}, Proposition 3.3). 
Note that the dual $\tilde{Y}_X^\vee$ has a point of multiplicity three, 
a special property already observed in \cite{EPW3}, Proposition 6.1, and that 
we have met in the proof of Proposition \ref{sextics}. 

\medskip What about the missing case $k=2$? And can we understand the relations between the
families of EPW sextics obtained from different values of $k$? 

To answer the latter question, we can use the construction 
of Gushel threefolds as degenerations of non Gushel Fano threefolds of degree ten 
\cite{gushel}. To be more specific, consider the projective cone over $G$, that we 
denote by $CG\subset\PP(\CC\oplus\wedge^2V_5)$. Let $p_0$ denote the 
vertex of this cone. Now we cut $CG$ by a general quadric $Q$, and a 
general linear 
space $\PP V_{10-k}$, of codimension $k+1$. We get a variety $Z$ of 
dimension $5-k$. There are two cases:
\begin{itemize}
\item $p_0\notin \PP V_{10-k}$: then $Z$ is isomorphic with the intersection 
$X$ of $G$ with the projection of $\PP V_{10-k}$ to $\PP(\wedge^2V_5)$, and 
a quadric $Q'$; 
\item $p_0\in \PP V_{10-k}$, that is, $\PP V_{10-k}$ is a cone over some 
$\PP V_{9-k}\subset \PP(\wedge^2V_5)$: then $Z$ is a double cover of 
$G\cap \PP V_{9-k}$, branched over its intersection $X$ with a quadric $Q'$. 
\end{itemize}
Obviously the second case is a degeneration of the first one. We call 
the corresponding $Z$ {\it Gushel varieties}. 

\begin{lemm}
If $Z$ is Gushel, the sextics $Y_Z$ and $Y_X$ are equal. 
\end{lemm}

\proof Take coordinates $(t,\omega)$ on $\CC\oplus\wedge^2V_5$ and write the 
equation of the quadric $Q$ as $Q(t,\omega)=t^2+2\ell(\omega)t+q(\omega)$. Note that
we can choose for the quadric $Q'$ defining $X$, the discriminant $Q'(\omega)=q(\omega)-\ell(\omega)^2$. 

Suppose that $Q+P_v$ defines a singular quadric in $\PP(\CC\oplus\wedge^2V_5)$. This means that 
we can find a point $(t_0,\omega_0)$ such that $Q((t_0,\omega_0),(t,\omega))+P_v(\omega_0,\omega)=0$
for any $(t,\omega)$. That is, we must have $t_0+\ell(\omega_0)=0$ and 
$$t_0\ell(\omega)+q(\omega_0,\omega)+P_v(\omega_0,\omega)=0$$
for any $\omega$. But this implies that $\omega_0$ defines a singular point of $Q'$. 
Hence $Y_Z\subset Y_X$, and since they are both sextics hypersurfaces, they are equal. \qed

\medskip Now consider the Gushel manifold $Z$ as a degeneration of a family of non Gushel manifolds. 
Suppose that $Z$ be defined by a quadric $Q(t,\omega)$ as above, 
and a linear space $\PP V_{10-k}$ through $p_0$, 
defined by the $k+1$ equations $h_0(\omega)=\cdots = h_k(\omega)=0$. Then we define $Z(\epsilon)$ 
by the same quadric, and the linear space $\PP V_{10-k}(\epsilon)$ with equations 
$$h_0(\omega)=\epsilon t, \qquad h_1(\omega)=\cdots = h_k(\omega)=0.$$
For $\epsilon\ne 0$, $\PP V_{10-k}(\epsilon)$ does not contain $p_0$. Hence $Z(\epsilon)$ is isomorphic
with the intersection $Z^*(\epsilon)$ of $G$ with the linear space $\PP V_{10-k}^*(\epsilon)$
of equations $h_1(\omega)=\cdots = h_k(\omega)=0$, and the quadric $Q(\epsilon^{-1}h_0(\omega),
\omega)=0$. 

\begin{lemm}
The sextic $Y_Z$ is a degeneration of the sextics $Y_{Z^*(\epsilon)}$.
\end{lemm}

\proof This is rather clear. The sextic $Y_{Z^*(\epsilon)}$ is defined by the condition that 
the quadric $zQ(\epsilon^{-1}h_0(\omega),\omega)+P_v(\omega)$ be singular on $\PP V_{10-k}^*(\epsilon)$, 
or equivalently, that the quadric $zQ(t,\omega)+P_v(\omega)$ be singular on $\PP V_{10-k}(\epsilon)$. 
Letting $\epsilon$ tend to zero, we get the sextic $Y_Z$ as a degeneration of the sextics 
$Y_{Z^*(\epsilon)}$. \qed

\medskip We can conclude inductively that for any $k\ge 0$, and any Fano manifold $X$ 
of degree ten and dimension $5-k$, the associated sextic $Y_X$ is a possibly degenerate
EPW sextic. 

\medskip\noindent {\it Remark}. For $k$ {\it odd}, the general quadric in $Y_X$, having 
corank one, is a cone over a smooth quadric of even dimension. Such a quadric has two rulings
by maximal linear subspaces, and this induces the double cover $\tilde{Y}_X\ra Y_X$, branched 
over the locus parametrizing quadrics of corank at least two. By the preceding construction,
this remark can be extended to the case where $k$ is even. The double covering $\tilde{Y}_X$
is endowed, by a deformation argument, with a symplectic structure, for $X$ general of any dimension.

\begin{prop}
For $X$ a general Fano threefold of degree ten, the associated sextic $Y_X$ is a general
EPW sextic. 
\end{prop}

\proof This follows from a dimension count. Remember that EPW sextics have 20 moduli. 
On the other hand, if $Y$ is an EPW sextic, and $X$ is a general Fano threefold such that $Y_X\simeq Y$, 
then we know from \cite{log2} that the singular locus of $Y$ is a smooth surface isomorphic 
with the Fano surface of conics in $X$ (more precisely, with the quotient of the minimal model 
of that surface, by a base point free involution). Moreover, Logachev's reconstruction theorem
(see the Appendix of \cite{dim})
implies that there is only a two dimensional family of Fano threefolds $X$ with the same 
Fano surface, and a fortiori with the same associated sextic $Y$. Since Fano threefolds
of degree ten have 22 moduli, this implies that $Y$ lives in a 20-dimensional family, 
hence must be a generic EPW sextic. \qed

\begin{coro}\label{2EPW}
For any $N=3,4,5$, and $X$ a general Fano manifold of degree ten and dimension $N$, 
the associated sextic $Y_X$ is a general EPW sextic.
\end{coro}

\proof This is a direct consequence of the previous degeneration
argument to Gushel type manifolds. \qed 

\medskip
The next obvious question to ask is: which are the Fano manifolds $X$ of degree ten and dimension $N$, 
whose associated sextic $Y_X$ is a given general EPW sextic $Y$? Denote by $m_N$ the dimension of 
the moduli space\footnote{We suppose implicitly that this moduli space does exist. We hope to
come back to this question in a future paper.} of Fano manifolds of degree ten and dimension $N$. 
An easy computation shows that 
$$m_2=19, \quad m_3=22, \quad m_4=24, \quad m_5=25.$$ 
The relative dimension $r_N$ of the map to the moduli space of EPW sextics is therefore given by
$$r_2=-1, \quad r_3=2, \quad r_4=4, \quad r_5=5.$$ 
For $N=3$, we have seen that the family $EPW^{-1}_N(Y)$ of Fano threefolds $X$ whose associated 
EPW sextic is isomorphic with  $Y$, is essentially the surface $S(Y)=Sing(Y^\vee)$. It is tempting 
to imagine that a similar phenomenon should hold for $N=4$ or $5$. 

\medskip\noindent {\it Question}. 
If $Y$ is a generic EPW sextic, is it true that 
$$EPW^{-1}_4(Y)\simeq Y^\vee-S(Y)\qquad \mathrm{and} \quad EPW^{-1}_5(Y)\simeq \PP^5-Y^\vee ?$$

\medskip
Indeed, for $N=4$, once we have a representation of $Y$ as $Y_Z$,  or equivalently, of $Y^\vee$ as 
$Y_Z^\vee$, the Pl\"ucker point, which belongs to  $Y^\vee-S(Y)$, is given a special role.  
We believe that specifying that point in $Y^\vee-S(Y)$ should be equivalent to specifying $Z$. 
The same phenomenon should hold for $N=5$, except that in that case the Pl\"ucker point does not belong
to $Y^\vee$. 

\subsection{O'Grady's double overs}
We can now prove that for $Z$ a general Fano fourfold of degree ten, our double cover $\tilde Y_Z^\vee$ 
of the general EPW sextic $Y_Z^\vee$ coincides with the double cover constructed by O'Grady 
(see \cite{EPW3}, \S 4). We denote the
latter by $\tilde Y_{Z,O}^\vee$. 
  
\begin{prop}\label{sameog}
The symplectic manifolds $\tilde Y_Z^\vee$ and $\tilde Y_{Z,O}^\vee$ are isomorphic. 
In particular $\tilde Y_Z^\vee$ is an irreducible symplectic manifold.
\end{prop}

\proof 
Recall that we denoted by $S_Z$ the singular locus of $Y_Z^\vee$. For a general $Z$ this is 
a smooth surface. Since $\tilde Y_{Z,O}^\vee$ is simply-connected, the \'etale double cover
$\tilde Y_{Z,O}^\vee-\tilde S_{Z,O}\ra Y_Z^\vee-S_Z$ (where $\tilde S_{Z,O}$ denotes the preimage
of $S_Z$) is the universal covering, and in particular $\pi_1(Y_Z^\vee-S_Z)=\ZZ_2$ (see
\cite{EPW4}, \S 3.2). Then since $\tilde Y_Z^\vee-\tilde S_Z\ra Y_Z^\vee-S_Z$ is also a 
non trivial  \'etale double cover, it lifts to an isomorphism between $\tilde Y_{Z,O}^\vee-\tilde S_{Z,O}$
and $\tilde Y_Z^\vee-\tilde S_Z$. So $\tilde Y_{Z,O}^\vee$ and $\tilde Y_Z^\vee$ are birational, 
and in particular $h^{2,0}(\tilde Y_Z^\vee)=h^{2,0}(\tilde Y_{Z,O}^\vee)=1$. 

Moreover, being birational, $\tilde Y_{Z,O}^\vee$ and $\tilde Y_Z^\vee$ are also deformation 
equivalent \cite[Theorem 4.6]{huy}. Therefore $\tilde Y_Z^\vee$ is, as $\tilde Y_{Z,O}^\vee$, 
a numerical $(K3)^{[2]}$, according to O'Grady's terminology. 
And we can conclude the proof by applying Theorem 1.1 in \cite{EPW3},
once we know that:


\begin{lemm}
The defining involution $\iota$ of the double covering 
$\tilde{Y}_Z^\vee\ra Y_Z^\vee$ is anti-symplectic.
\end{lemm}

\proof
Let $c$ be a general conic in $Z$, and $V_4$ the corresponding hyperplane 
in $V_5$. Let $y$ be the image of $[c]$ in $\tilde{Y}_Z^\vee$. Recall that 
the quartic surface $S_c=G(2,V_4)\cap Q\cap H$ has normal bundle $N_{S_c/Z}
\simeq T^\vee_{S_c}$. Moreover the normal 
sequence to the triple $(c\subset S_c\subset Z)$, gives rise to the
exact sequence
$$0\ra H^0(N_{c/S_c})\ra H^0(N_{c/Z})\ra H^0(N_{S_c/Z|c})\ra 0,$$
which must be interpreted as the tangent sequence of the map 
$F_g(Z)\ra\tilde Y_Z^\vee$. In particular, this yields 
a natural identification 
$$T_y\tilde{Y}_Z^\vee\simeq H^0(N_{S_c/Z|c})\simeq H^0(T_c^\vee)\simeq V_4^\vee,$$
and the symplectic form $\tilde\Phi_\sigma$ can be defined at $y$ as the 
composition 
$$\wedge^2T_y\tilde{Y}_Z^\vee\simeq\wedge^2H^0(T_c^\vee)\ra 
H^0(\wedge^2T_c^\vee)=H^0(\cO_c(2))\stackrel{\kappa_c}{\ra} \CC.$$
Here, recall that the quadratic form $\kappa_c$ is induced by the twisted conormal 
sequence
$$0\ra\omega_c\ra N^\vee_{c/Q_{c,V_4}}(2)\ra \cO_c(2)=\cO_Z(1)_{|c}\ra 0,$$
where $Q_{c,V_4}=G(2,V_4)\cap H$. 

Take another conic $c'$ in $Z$ such that the  image of $[c']$  
in $\tilde{Y}_Z^\vee$ be the point $y'=\iota(y)$. Recall that this means 
that the planes $\langle c\rangle$ and $\langle c'\rangle$ spanned by the 
two conics are contained 
in the same (singular) quadric  $Q_{c,V_4}=Q_{c',V_4}$ of the pencil of quadrics 
we have in $\PP(\wedge^2V_4)\cap H$, but do not belong to the same ruling. 
In particular the two planes $\langle c\rangle$ and $\langle c'\rangle$
meet along a line, and the three-plane $\langle c,c'\rangle$ cuts 
$Z$ along the degenerate elliptic curve 
$$e =c\cup c'=G(2,V_4)\cap Q\cap \langle c,c'\rangle.$$
Now, if we take an element of $\wedge^2V_4^\vee$ and apply 
$\tilde\Phi_{\sigma,y}+\tilde\Phi_{\sigma,y'}$, we first get an element 
of $H^0(\cO_c(2))\oplus H^0(\cO_{c'}(2))$ which defines a global section 
of $H^0(\cO_Z(1)_{|e})$. Then we apply the coboundary maps defined 
by the twisted conormal sequences of $c$ and $c'$. But the analogous
conormal sequence for $e$, 
$$0\ra\omega_e\ra N^\vee_{e/Q_e}(2)\ra 
\cO_Z(1)_{|e}\ra 0,$$
splits since $e$ is a complete intersection curve in $Q_e=Q_{c,V_4}=Q_{c',V_4}$. 
This implies that the associated coboundary operator is zero, which means that 
$\tilde\Phi_{\sigma,y}+\tilde\Phi_{\sigma,y'}$ is zero. Otherwise stated, 
$$\tilde\Phi_\sigma+\iota^*\tilde\Phi_\sigma=0,$$
or else, $\iota$ is anti-symplectic.
 \qed

\section{Two integrable systems}

\subsection{Fano threefolds contained in $Z$}

Let $W=G\cap Q\cap H'\cap H''$ be a general Fano threefold of degree ten 
contained in the general  Fano fourfold $Z=G\cap Q\cap H$ as a hyperplane 
section. The linear systems $I_2(Z)$ and $I_2(W)$ of quadrics containing 
them are naturally identified. As before, we associate to $W$ the hypersurface
$\hat{Y}_W^\vee\subset I^\vee$ and its singular locus, the surface $S_W$.

\begin{lemm}
The surface $S_W$ is contained in $Y_Z^\vee$.
\end{lemm}

\proof By definition, a point $h$ in $S_W$ is a hyperplane in $I$ 
such that for some hyperplane $V_4$ of $V_5$, quadrics in $h$ restrict 
to a corank two quadric in $\PP(\wedge^2V_4)\cap H'\cap H''$. This 
is a hyperplane in $\PP(\wedge^2V_4)\cap H$, and it follows that quadrics
in $h$ cut $\PP(\wedge^2V_4)\cap H$ along a quadric with some corank 
two hyperplane section. But this is possible only if  that quadric is
itself singular. This precisely means that $h$ defines a point of   
$Y_Z^\vee$.\qed

\medskip
As explained in \cite{log2}, and as we already mentioned, 
the surface $S_W$ is closely related to 
conics in $W$. Indeed, the Hilbert scheme $F_g(W)$ parametrizing 
conics in $W$ is a smooth surface, containing a unique $\rho$-conic,
and a line of $\sigma$-conics. This line is an exceptional curve which 
may be contracted. The resulting surface $F_m(W)$ is then endowed 
with a fixed-point free involution whose quotient is precisely $S_W=
F_\iota(W)$. It follows that $F_m(W)$ can be seen as the pull-back 
$\tilde{S}_W$ of  $S_W$ inside $\tilde{Y}_Z^\vee$, and that we have 
a commutative diagram
$$\begin{array}{ccccc}
F_g(W) & \ra & F_m(W)=\tilde{S}_W & \ra & F_\iota(W)=S_W \\
\downarrow & & \downarrow & & \downarrow \\
F_g(Z) & \ra & \tilde{Y}_Z^\vee & \ra & Y_Z^\vee, 
\end{array}$$
where the vertical maps are injections. 

\begin{prop}\label{lag1}
The surface $\tilde{S}_W$ is a Lagrangian subvariety of $\tilde{Y}_Z^\vee$. 
\end{prop}

\proof We just need to prove that $F_g(W)$ is isotropic with respect to
the two-form $\phi_\sigma$ on $F_g(Z)$. Otherwise said, for a general 
conic $c$ in $W$, we must check that $\phi_\sigma$ vanishes on 
the subspace $T_{[c]}F_g(W)=H^0(N_{c/W})$ of $T_{[c]}F_g(Z)=H^0(N_{c/Z})$. 
Since $W$ is a hyperplane section of $Z$, the conormal sequence of 
the triple $(c,W,Z)$ is just
$$0\ra N_{c/W}\ra N_{c/Z}\ra \cO_c(2)\ra 0.$$
In particular $N_{c/W}$ has degree zero (in fact $N_{c/W}$ is trivial for a
general conic, but we will not need that -- see \cite{dim} for more details).
We have a commutative diagram
$$\begin{array}{ccccc}
\wedge^2H^0(N_{c/W}) & \hookrightarrow & \wedge^2H^0(N_{c/Z}) & & \\
\downarrow & & \downarrow & & \\
H^0(\wedge^2N_{c/W}) & \ra & H^0(\wedge^2N_{c/Z}) &  & \\
|| & & || & & \\
H^0(\cO_c) & \hookrightarrow & H^0(N_{c/Z}^\vee(2)) & \ra & H^1(\omega_c)=\CC ,
\end{array}$$
and we need to prove that the composition $\wedge^2H^0(N_{c/W})\ra
H^1(\omega_c)$ is zero. But recall that the map $H^0(N_{c/Z}^\vee(2)) 
\ra H^1(\omega_c)$ was induced by the twisted conormal sequence of 
the triple $(c,Z,X=G\cap H)$. This sequence fits with the conormal
sequence of the triple $(c,W,Y=G\cap H'\cap H'')$ into the following
commutative diagram:
$$\begin{array}{ccccccccc}
 &&&&0&&0&&\\
 &&&&\downarrow &&\downarrow &&\\
 &&&&\cO_c&=&\cO_c&&\\
 &&&&\downarrow &&\downarrow &&\\
0&\ra &\omega_c&\ra & N_{c/X}^\vee(2)&\ra & N_{c/Z}^\vee(2)&\ra & 0 \\
&&||&&\downarrow &&\downarrow &&\\
0&\ra &\omega_c&\ra & N_{c/Y}^\vee(2)&\ra & N_{c/W}^\vee(2)&\ra & 0 \\
 &&&&\downarrow &&\downarrow &&\\
 &&&&0&&0&&
\end{array}$$
The map $H^0(N_{c/Z}^\vee(2)) \ra H^1(\omega_c)$ is the coboundary map 
of the middle exact sequence. But this diagram shows that the sequence
is split over the factor $\cO_c=\wedge^2N_{c/W}$ of $N_{c/Z}^\vee(2)$.
Therefore the coboundary map vanishes on $H^0(\cO_c)\subset 
H^0(N_{c/Z}^\vee(2))$, and our claim follows.\qed

\medskip
We are thus in the situation where we can use the results of Donagi 
and Markman about deformations of a Lagrangian subvariety $S$ of 
a symplectic variety $Y$ \cite{dm}: over the Hilbert scheme $\cB$ parametrizing
smooth deformations of $S$ in $Y$ (which are non-obstructed), there exists
an integrable system, otherwise said a Lagrangian fibration, 
whose Liouville tori are the Albanese varieties $Alb(S)$.

\smallskip
In our setting, note that the Abel-Jacobi mapping 
$AJ : F_g(W)\ra J(W)$ factorizes through
$F_m(W)=\tilde{S}_W$ and induces an isomorphism (see \cite{log2,dim})
$$alb(AJ) : \;\; Alb(\tilde{S}_W)\simeq J(W).$$

Denote by $U_Z\subset\PP V_9^\vee$ the open subset parametrizing 
smooth hyperplane section $W$ of $Z$. 
We deduce the following statement:

\begin{theo}\label{int1}
For a general Fano fourfold $Z=G\cap H\cap Q$ of degree ten,
the set $U_Z$ parametrizing smooth Fano threefolds $W=G\cap H'\cap H''\cap Q$
in $Z$, is contained in the base $\cB$ of an integrable system, in such a 
way that over $U_Z$  
the Liouville tori are the intermediate Jacobians $J(W)$.
\end{theo}

An interesting point here is that $U_Z$ has dimension eight, while $\cB$ is ten 
dimensional. In particular, the deformations of $S=\tilde{S}_W$ in $Y=\tilde{Y}_Z^\vee$
are not all obtained by deforming $W$ in $Z$. This is certainly related to the fact
that the representation $Y=\tilde{Y}_Z^\vee$ does not defined $Z$ uniquely, as we 
have stressed in section 4.5. Deforming $Z$ without changing $Y$, and taking 
hyperplane sections, we should get more deformations of $S$. 

To be more specific, we can observe that the representation $Y=\tilde{Y}_Z^\vee$ 
gives a special role to the two Pl\"ucker points (the two preimages of the  
Pl\"ucker point in $Y_Z^\vee$), and that the surfaces $\tilde{S}_W$ always contain
these points (since $W$ always contain conics of type $\rho$ or $\sigma$). 
Therefore, deforming  $W$ in $Z$ should be equivalent to deforming $S=\tilde{S}_W$ in 
$Y=\tilde{Y}_Z^\vee$ with the Pl\"ucker points fixed. 

\smallskip

\subsection{Fano fivefolds containing $Z$}

Now we consider the moduli stack $\cB$ parametrizing smooth 
fivefolds $X=G\cap Q$ containing a fixed 
fourfold $Z=G\cap Q\cap H$ as a hyperplane section. 
By \cite{beau}, the tangent space to $\cB$ at the point defined
by the fivefold $Z$ can be identified with $H^1(X,TX(-1))$. 

\begin{lemm}
The Zariski tangent space $H^1(X,TX(-1))$ to $\cB$ at $[X]$ is 
naturally isomorphic with $H^2(X,\Omega^3_X)$. Its dimension is ten.
\end{lemm}

\proof 
Since $\omega_X=\cO_X(-3)$, we have $H^2(X,\Omega^3_X)=H^2(X,\wedge^2TX(-3))$. 
The normal exact sequence of the inclusion $X\subset G$ induces the 
exact sequence
$$0\ra\wedge^2TX\ra\wedge^2TG_{|X}\ra TX(2)\ra 0.$$
By Bott's theorem $\wedge^2TG(-3)$ is acyclic, and $\wedge^2TG(-5)=\Omega^4_G$
has non zero cohomology only in degree four. Therefore $\wedge^2TG(-3)_{|X}$
has non zero cohomology only in degree three. This implies that 
$H^1(X,TX(-1))\simeq H^2(X,\wedge^2TX(-3)$, as claimed. \qed

\medskip
Now consider the EPW sextic $Y_X$ and its singular locus $S_X$, which 
is a smooth surface. 

\begin{prop}
The surface $S_X$ is contained in $\hat{Y}_Z^\vee$.
\end{prop}

\proof Recall that $S_X$ parametrizes pairs $(h,V_4)$ made of 
hyperplanes $h$ in $I_2(X)$, and hyperplanes
$V_4\subset V_5$, such that the pencil of quadrics on $\PP(\wedge^2V_4)$
obtained by restricting $I_2(X)$, contains a quadric of rank four 
whose preimage in $I_2(X)$ is precisely $h$.
Cutting with the hyperplane $H$ defining $Z$, we remain with a quadric
of rank at most four, which implies that the point $(h,V_4)$ belongs 
to $\hat{Y}_Z^\vee$. \qed

\medskip
Now we lift this surface to $\tilde{Y}_Z^\vee$. We get the diagram:

$$\begin{array}{ccc}
\tilde{S}_X & \hookrightarrow & \tilde{Y}_Z^\vee \\
\downarrow &&\downarrow \\
S_X & \hookrightarrow & \hat{Y}_Z^\vee .
\end{array}$$

\smallskip

\begin{prop}\label{lag2}
The surface $\tilde{S}_X$ is a Lagrangian subvariety of $\tilde{Y}_Z^\vee$.
\end{prop}

\proof A point $(h,V_4)\in S_X$ defines a corank two quadric in a 
$\PP(\wedge^2V_4)$, and such a quadric contains two pencils of 
three-planes. Each of these three-planes will cut $X$ along quadratic 
surface. 

A general quadratic surface $\Sigma$ in $G$ is given by two transverse planes 
$V_2$ and $V'_2$, as the image of the obvious map $\PP(V_2)\times \PP(V'_2)
\hookrightarrow G(2,V_2\oplus V'_2)\subset G$. The corresponding parameter 
space is and open subset of $Sym^2G$ and has dimension $12$. The normal 
bundle of $\Sigma$ in $G$ is easily seen to decompose as
$$N_{\Sigma/G}=\cO_\Sigma(1,1)^{\oplus 2}\oplus 
\cO_\Sigma(1,0)\oplus \cO_\Sigma(0,1).$$
If $\Sigma\subset X$, its normal bundle is the kernel of the 
induced exact sequence 
$$0\ra N_{\Sigma/X}\ra N_{\Sigma/G}\ra \cO_\Sigma(2,2)\ra 0.$$
Generically $h^1(N_{\Sigma/X})=0$ and $h^0(N_{\Sigma/X})=12-3\times 3=3$, 
so that there is a smooth three-dimensional family of quadratic surfaces
in $X$. There is a natural map from this family to $F_g(Z)$, defined by 
cutting a quadratic surface $\Sigma$ with the hyperplane $H$ spanned by $Z$, 
to get a conic $c=\Sigma\cap H$. We are then reduced to showing 
that the image of the restriction map 
$$H^0(N_{\Sigma/X})\ra H^0(N_{\Sigma/X|c})=H^0(N_{c/Z})$$
is isotropic with respect to the two-form $\phi_\sigma$.

But this is easy: recall that if $Y=G\cap H$, so that $Z=Y\cap Q$, 
the two-form $\phi_\sigma$ was defined with the help of the 
normal exact sequence of the triple $(c,Z,Y)$. But this is the restriction 
to $c$ of the normal exact sequence of the triple $(\Sigma,X,G)$, 
which reads, after dualizing and twisting, 
$$0\ra\cO_\Sigma(-1,-1)\ra N_{\Sigma/G}^\vee(1)\ra N_{\Sigma/X}^\vee(1)\ra 0.$$
Otherwise said, there is a commutative diagram
$$\begin{array}{ccccc}
\wedge^2H^0(N_{c/Z}) & \ra & H^0(\wedge^2N_{c/Z})=H^0(N_{c/Z}^\vee(1)) & 
\ra & H^1(\omega_c)=\CC \\
\uparrow & & \uparrow & & \uparrow \\
\wedge^2H^0(N_{\Sigma/X}) & \ra & H^0(\wedge^2N_{\Sigma/X})=
H^0(N_{\Sigma/X}^\vee(1)) & 
\ra & H^1(\cO_\Sigma(-1,-1)). 
\end{array}$$
The first line defines $\phi_\sigma$, and the last line, its restriction 
to $H^0(N_{\Sigma/X})$. Since $H^1(\cO_\Sigma(-1,-1))=0$, this restriction 
vanishes, and we are done. 
\qed

\begin{theo}\label{int2}
For a general Fano fourfold $Z=G\cap H\cap Q$ of degree ten,
the moduli stack $\cB$ parametrizing smooth Fano fivefolds $X=G\cap Q$
containing $Z$ is the base of an integrable system 
whose Liouville tori are the intermediate Jacobians $J(X)$.
\end{theo}

\proof This can be proved as in \cite{m1} for K3-Fano flags, or as in
\cite{m2} for cubic fivefolds containing a given cubic fourfold. Let us 
briefly recall the argument, which goes back to \cite{dm}, with the necessary 
(minor) modifications. 

A first observation is that the normal exact sequence of the pair $(Z,X)$ 
induces isomorphisms
$$H^1(\Omega^4_X(Z))\simeq H^1(\Omega^4_X(Z)_{|Z})\simeq H^1(\Omega^3_Z)\simeq\CC.$$
(For the first two isomorphisms, there are some easy vanishing to verify. For the last
one see Lemma \ref{hodge}.) One then checks that tensoring with a generator $\omega_Z$ of $H^1(\Omega^4_X(Z))$
defines an isomorphism
$$H^1(TX(-Z))\simeq H^2(\Omega^3_X).$$
The left hand side is to be interpreted as the tangent space to $\cB$ at the point
defined by $Z$. The right hand side is the fiber of the Hodge bundle $\cH^{3,2}(\cX/\cB)$. 
The dual vector bundle $\cE$ on $\cB$ is thus endowed with a natural symplectic form, 
and one must check that this form descends to the intermediate Jacobian bundle, its 
quotient by the locally constant bundle of integral forms. For this, one has to 
normalize the isomorphism $\cH^{3,2}(\cX/\cB)\simeq \Omega^1_\cB$ by requiring 
that over $Z$, it is defined by a generator $\omega_Z$ of $H^1(\Omega^4_X(Z))$
restricting to a fixed generator of $H^1(\Omega^3_Z)$. Then the proof
\cite{m2}, Theorem 2.3, applies verbatim. \qed

\medskip 
It is probably possible to deduce Theorem \ref{int2}  directly from 
Proposition \ref{lag2}, as we deduced Theorem \ref{int1} from Proposition \ref{lag1}.
Indeed the general results of Donagi-Markman imply that one can define over $\cB$ 
an integrable system whose fiber over $X$ is the Albanese variety $Alb(\tilde{S}_X)$. 

On the other hand, consider the Hilbert scheme 
$F_{qs}(X)$ parametrizing quadratic surfaces in $X$. Once a point is chosen 
in this scheme, the Abel-Jacobi mapping gives a morphism 
$$AJ : F_{qs}(X)\ra J(X).$$
By the previous arguments $F_{qs}(X)$ is a $\PP^1$-bundle over $\tilde{S}_X$,
and since every map from $\PP^1$ to a complex torus is constant, 
we get an induced morphism $AJ : \tilde{S}_X\ra J(X)$. Hence, for the 
Albanese variety, a morphism
$$alb(AJ) : \;\; Alb(\tilde{S}_X)\ra J(X).$$
Very probably, this morphism should be an isomorphism. But this seems 
technically much more difficult to check than to prove Theorem \ref{int2}
as we did above. 


\section*{Appendix : Proof of Lemma \ref{double-sing-tau}}

We can choose a basis of $V_5$ such that $\ell$ be the double line
defined as the intersection of $G(2,V_5)$ with the plane $P=\langle  
v_1\wedge v_2,  v_1\wedge v_3, v_2\wedge v_3+v_1\wedge v_4\rangle$.
Around $v_1\wedge v_3$ we have affine coordinates on $G(2,V_5)$ such that 
a plane transverse to $\langle v_2,v_4,v_5\rangle$ has a basis of the form
\begin{eqnarray*}
 w_1 &= v_1+ z_2v_2+z_4v_4+z_5v_5, \\
 w_3 &= v_3+ t_2v_2+t_4v_4+t_5v_5. 
\end{eqnarray*}
In these coordinates, $\cI_{\ell,G}$ is generated by $z_4,z_5,t_4-z_2,t_5$ and
$t_4^2$. We check that an element $\psi\in Hom(\cI_{\ell,G},\cO_\ell)$ is of 
the following form:
\begin{eqnarray*}
t_4^2 &\mapsto & \psi_1+\psi_2t_2+\psi_3t_2^2+\psi_4t_4+\psi_5t_2t_4, \\
t_5 &\mapsto & \psi_6+\psi_7t_2+\psi_8t_4, \\
t_4-z_2 &\mapsto & \psi_9+\psi_{10}t_2+\psi_{11}t_4, \\
z_5 &\mapsto & \psi_{12}+\psi_7t_4, \\
z_4 &\mapsto & \psi_{13}+\psi_3t_2+(\psi_5-\psi_{10})t_4.
\end{eqnarray*}
Now, locally around $e_1\wedge e_3$ the ideal sheaf $\cI_{Z,G}$ is 
generated by $h/p_{13}$ and $q/p_{13}^2$, where $h$ and $q$ are the equations
of the hyperplane $H$ and of the quadric $Q$, expressed in terms of Pl\"ucker 
coordinates. Write
\begin{eqnarray*}
h =\sum_{i<j}h_{ij}p_{ij}, \qquad
q = \sum_{i<j,k<\ell}q_{ij,kl}p_{ij}p_{kl}.
\end{eqnarray*}
For $Z$ to contain $\ell$ we first need that $H$ contains $P$, that is,
$$h_{12}=h_{13}=h_{14}+h_{23}=0,$$
and that the equation of $Q$ restricted to the plane 
$P$ reduces to that of the double line $\ell$, which gives
\begin{eqnarray*}
 &q_{12,12}=q_{12,13}=q_{13,13}=0, \\
 &q_{12,14}+q_{12,23}=q_{13,14}+q_{13,23}=0.
\end{eqnarray*}
Using these relations, we must express $h/p_{13}$ and $q/p_{13}^2$
in terms of our preferred local generators of $\cI_{\ell, G}$ and deduce
their images by the morphism $\psi$. We find that 
\begin{eqnarray*}
h/p_{13} &\mapsto & A+Bt_2+Ct_4, \\
q/p_{13}^2 &\mapsto & D+Et_2+Ft_2^2+Gt_4+Ht_2t_4
\end{eqnarray*}
where the quantities $A,B,C,D,E,F,G,H$ are given by the following formulas:
\begin{eqnarray*}
A &= &h_{15}\psi_6+h_{14}\psi_9+h_{24}\psi_1-
h_{34}\psi_{13}-h_{35}\psi_{12}, \\
B &= &h_{15}\psi_7+h_{14}\psi_{10}+h_{24}(\psi_2-\psi_{13})-
h_{34}\psi_3-h_{25}\psi_{12}, \\ 
C &= &h_{15}\psi_3+h_{14}\psi_{11}+h_{24}(\psi_4-\psi_9)-
h_{34}(\psi_5-\psi_{10}) \\
 & & \hspace*{6cm}-h_{35}\psi_7+h_{25}\psi_6-h_{45}\psi_{12}, \\
D &= &b_{15}\psi_6+(b_{24}+g)\psi_1-b_{34}\psi_{13}-
b_{35}\psi_{12}+e\psi_{9}, \\
E &= &a_{15}\psi_6+b_{15}\psi_7+a_{24}\psi_1+b_{24}(\psi_2-\psi_{13})
-(a_{35}+b_{25})\psi_{12} \\ 
 & & \hspace*{4cm}-a_{34}\psi_{13}-b_{34}\psi_{3}
+d\psi_{9}+e\psi_{10}+g\psi_{2}, \\
F &= &a_{15}\psi_7+a_{24}(\psi_2-\psi_{13})-a_{25}\psi_{12}
+(g-a_{34})\psi_3+d\psi_{10}, \\
G &= &b_{15}\psi_8+(b_{25}+c_{15})\psi_6+c_{24}\psi_1+b_{24}(\psi_4-\psi_{9})
-b_{34}(\psi_{5}-\psi_{10}) \\
 & & \hspace*{3cm}-c_{34}\psi_{13}-b_{35}\psi_{7}-c_{35}\psi_{12}
+f\psi_{9}+e\psi_{11}+g\psi_{4}, \\
H &= &a_{15}\psi_8+c_{15}\psi_7+a_{24}(\psi_4-\psi_{9})
+c_{24}(\psi_2-\psi_{13})+a_{25}\psi_6-c_{25}\psi_{12} \\
 & & \hspace*{2cm}-a_{34}(\psi_{5}-\psi_{10})-c_{34}\psi_3-a_{35}\psi_{7}
+g\psi_{5}+f\psi_{10}+d\psi_{11}.
\end{eqnarray*}
In these formulas we have set $a_{ij}=q_{12,ij}$, $b_{ij}=q_{13,ij}$,
$c_{ij}=2q_{14,ij}=2q_{23,ij}$, $d=q_{12,14}=-q_{12,23}$, $e=q_{13,14}=
-q_{13,23}$, $f=-2q_{23,23}$ and $g=q_{14,14}+q_{23,23}$.

So the rank of $\phi$ is equal to the rank of the $13\times 8$ matrix
$$\small{
\begin{pmatrix}
b_{24}+g & a_{24} & 0  & c_{24} & 0  & h_{24} & 0  & 0 \\
0 & b_{24}+g & a_{24} & 0  & c_{24} & 0  & h_{24} & 0  \\
0 & -b_{34} & g-a_{34} & 0  & -c_{34} & 0  & -h_{34} & 0  \\
0 & 0 & 0 & b_{24}+g & a_{24} & 0 & 0  & h_{24} \\
0 & 0 & 0 & -b_{34} & g-a_{34} & 0  & 0  & -h_{34} \\
b_{15} & a_{15} & 0 & c_{15}+b_{15} & -a_{25} & h_{15} & 0 & h_{25} \\
0 & b_{15} & a_{15} & -b_{35} & c_{15}-a_{35} & 0 & h_{15} & -h_{35} \\
0 & 0 & 0 & b_{15} & a_{15} & 0 & 0 & h_{15} \\
e & d & 0 & f-b_{24} & -a_{24} & h_{14} & 0  & -h_{24} \\
0 & e & d & b_{34} & f+a_{34} & 0 & h_{14} & h_{34} \\
0 & 0 & 0 & e & d & 0 & 0 & h_{14} \\
-b_{35} & -b_{25}-a_{35} & -a_{25} &-c_{35} &-c_{25} & -h_{35}
 & -h_{25}  & -h_{45} \\
 -b_{34} & -b_{24}-a_{34} & -a_{24} &-c_{34} &-c_{24} & -h_{34}
 & -h_{24}  & 0 
\end{pmatrix}}$$
We need to show that this matrix has full rank outside a locus of codimension
at least three. For this we let $d_{24}=b_{24}+g$, $d_{34}=a_{34}+f$, $k=
-a_{34}-b_{24}$. After permuting lines and columns we get the following 
matrix $M$:
$$\small{
\begin{pmatrix}
d_{24} & a_{24} & 0  & h_{24} & 0 & c_{24} & 0    & 0 \\
0 & d_{24} & a_{24} & 0  & h_{24} & 0  & c_{24} & 0  \\
0 & 0 & 0 & 0 & 0 & d_{24} & a_{24}  & h_{24} \\
b_{15} & a_{15} & 0 & h_{15} & 0 &c_{15}+b_{15} & -a_{25} &  h_{25} \\
0 & b_{15} & a_{15} & 0 & h_{15} & -b_{35} & c_{15}-a_{35} & -h_{35} \\
0 & 0 & 0 & 0 & 0 &b_{15} & a_{15} &  h_{15} \\
e & d & 0 & h_{14} & 0  & d_{24}+k & -a_{24} &  -h_{24} \\
0 & e & d & 0 & h_{14} & b_{34} & d_{34} & h_{34} \\
0 & 0 & 0 & 0 & 0 &e & d &  h_{14} \\
0 & 0 & 0 & 0  & 0 & -b_{34} & d_{24}+k   & -h_{34} \\
-b_{35} & -b_{25}-a_{35} & -a_{25} &-h_{35}
 & -h_{25}  &-c_{35} &-c_{25} &  -h_{45} \\
 -b_{34} & k & -a_{24} &-h_{34}
 & -h_{24} & -c_{34} &-c_{24} &  0 \\
0 & -b_{34} & d_{24}+k & 0  & -h_{34} &0  & -c_{34} &  0  
\end{pmatrix}}$$

Observe the role of the matrix 
$$m=\begin{pmatrix}
d_{24} & a_{24} & h_{24} \\
b_{15} & a_{15} & h_{15} \\
e & d  & h_{14} 
\end{pmatrix}$$
Its rank is at least two in codimension three. If it  is equal to three, then clearly 
$\phi$ has full rank. So we may suppose that the rank of $m$ is equal to two, which occurs
in codimension one. Let $A,B,C$ denote the three-dimensional spaces corresponding 
to columns $124$, $235$ and $678$ respectively. Observing the three first lines
of the matrix, we see that they can be written as $p(v_1)+v_1', q(v_1)+v_1'', v_1$, 
where $v_1,v_1',v_1''$ belong to $C$ and $p:C\ra A$, $q:C\ra B$ are isomorphisms. 
Moreover the same is true for the two next groups of three lines, for some vectors
$v_2,v_3,..$ in $C$. Our hypothesis on $m$ is that the span of $v_1,v_2,v_3$ is 
two-dimensional. So there is a  relation $\alpha_1v_1+\alpha_2v_2+\alpha_3v_3=0$, and combining 
our lines accordingly we get the vectors $\alpha_1v'_1+\alpha_2v'_2+\alpha_3v'_3$
and $\alpha_1v''_1+\alpha_2v''_2+\alpha_3v''_3$, which belong to $C$. Since the 
tenth line of the matrix $M$ is also a vector of $C$, it is easy to conclude that 
in codimension at least three, $C$ is contained is the span of the lines of $M$. 

Then we can focus on the first five columns and forget the other ones. We know that 
the first nine lines of $M$ span a space of the form $p(L)+q(L)$ for some plane 
$L$ in $C$. If $p(L)$ and $q(L)$ meet, this can be only inside $A\cap B$, which is 
one dimensional. This is easy to exclude in codimension two. Then $p(L)+q(L)$ has 
codimension one, and to ensure that the matrix has full rank, it is enough to check
that the  the last three lines contribute, that is, they are not contained in 
$p(L)+q(L)$. Since the entries of these lines do not appear in the remaining of the
matrix, except for $a_{24}$ and $d_{24}$, this is also easy. \qed

\end{document}